\documentclass[a4paper,fleqn]{article}

\usepackage[utf8]{inputenc}
\usepackage[english]{babel}

\usepackage{graphicx}
\graphicspath{{fig/}{../fig/}}

\usepackage{csquotes}
\usepackage{xspace}

\usepackage{amsfonts,amsmath,amssymb}
\usepackage[usenames,dvipsnames]{xcolor}

\usepackage[space]{grffile}
\usepackage{latexsym}
\usepackage{textcomp}
\usepackage{longtable}
\usepackage{tabulary}
\usepackage{booktabs,array,multirow}

\usepackage{url}

\AtBeginDocument{\DeclareGraphicsExtensions{.pdf,.PDF,.eps,.EPS,.png,.PNG,.tif,.TIF,.jpg,.JPG,.jpeg,.JPEG}}

\usepackage{siunitx}
\usepackage{placeins}
\usepackage{subcaption}
\usepackage{subfiles}
\usepackage{wrapfig}
\usepackage{diagbox}
\usepackage{todonotes}

\DeclareMathOperator{\re}{Re}
\DeclareMathOperator{\im}{Im}
\renewcommand{\d}{\mathrm{d}}
\renewcommand{\i}{\mathrm{i}}
\newcommand{\e}{\mathrm{e}}

\definecolor{darkred}{RGB}{150,0,0}
\definecolor{darkblue}{RGB}{0,0,100}
\definecolor{darkgreen}{RGB}{0,80,0}

\usepackage[top=1in, bottom=1.25in, left=1.25in, right=1.25in]{geometry}

\newcommand{\beginsupplement}{%
	\setcounter{table}{0}
	\renewcommand{\thetable}{S\arabic{table}}%
	\setcounter{figure}{0}
	\renewcommand{\thefigure}{S\arabic{figure}}%
}

\begin{document}

\title{A Delay Equation Model for the Atlantic Multidecadal Oscillation}

\author{
	Swinda K.J. Falkena$^{1,2}$, Courtney Quinn$^{3,4}$, Jan Sieber$^4$, 
	and Henk A. Dijkstra$^{2,5}$}
\date{}

\maketitle

\noindent{$^{1}$Department of Mathematics and Statistics, University of Reading, Reading, UK\\
	$^{2}$Institute for Marine and Atmospheric Research Utrecht, Department of Physics, Utrecht University, Utrecht, The Netherlands\\
	$^{3}$CSIRO Oceans and Atmosphere, Hobart, TAS, AU\\
	$^{4}$Department of Mathematics, University of Exeter, Exeter, UK\\
	$^{5}$Centre for Complex Systems Studies, Faculty of Science, Utrecht University, Utrecht, The Netherlands\\
}

{\centering
	\textbf{Corresponding author}: {s.k.j.falkena@pgr.reading.ac.uk}\\[10mm]
}

\begin{abstract}
     A new technique to derive delay models from systems of partial differential equations, based on the Mori-Zwanzig formalism, is used to derive a delay difference equation model for the Atlantic Multidecadal Oscillation. The Mori-Zwanzig formalism gives a rewriting of the original system of equations which contains a memory term. This memory term can be related to a delay term in a resulting delay equation. Here the technique is applied to an idealized, but  spatially extended, model of the Atlantic Multidecadal Oscillation. The resulting  delay difference model is of a different type than the delay differential model which has been used to describe the El Ni\~no- Southern Oscillation. In addition to this model, which can also be obtained by integration along characteristics, error terms for a smoothing approximation of the model have been derived from the Mori-Zwanzig formalism.  Our new method of deriving delay models from spatially extended models has a large potential to use delay models to study a range of climate variability phenomena. 
\end{abstract}

\maketitle

\clearpage

\section{Introduction}
\label{sec:intro}

To better understand climate variability and climate change often conceptual climate models are used. These models  capture  the dominant physical processes behind climate phenomena, allowing for an improved understanding. Delay equation models form one class of conceptual climate models. These type of models are infinite dimensional, but often depend only on a few variables and parameters. This means that they can potentially describe more complex behaviour compared to ordinary differential equation (ODE) models, while still being easier to study than multi-dimensional partial differential equation (PDE) models.

Delay models have already been used to describe certain climate phenomena, particularly for the El Ni\~no-Southern Oscillation (ENSO) and Earth's Energy Balance  \cite{Keane2017}. Recently a new method of deriving delay equation models has been proposed \cite{Falkena2019}, allowing for a potential extension of the use of delay models to study other climate phenomena. This method of deriving delay models is based on the Mori-Zwanzig (MZ) formalism, which allows for the reduction of  high-dimensional systems to reduced-order models \cite{Chorin2002}. These reduced-order models are simpler to study, while still describing the physical processes present in the original high-dimensional model. So far the method in \cite{Falkena2019} has only been applied to  a PDE model of ENSO, for which ad-hoc delay models were already proposed \cite{Suarez1988}. Here we apply the MZ formalism to a PDE model of the Atlantic Multidecadal Oscillation (AMO), to investigate whether this phenomena can be described by a delay model as well.

\begin{figure}[h]
	\centering
	\includegraphics[width=.8\textwidth]{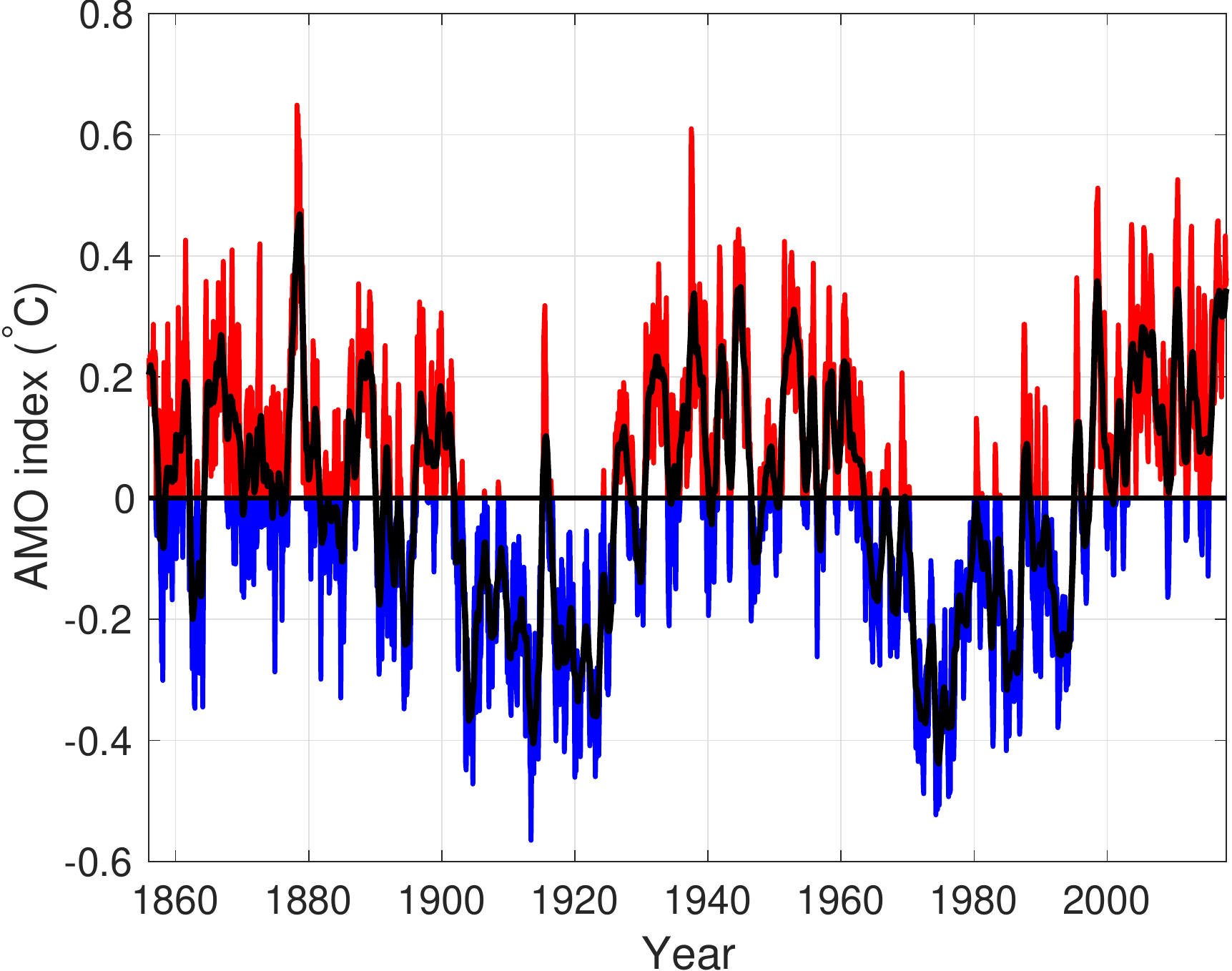}
	\caption{\label{fig:amoosc} The AMO index for the last 160 years. The index is computed as the deviations of the area-weighted average SST over the North Atlantic. In black the 12-monthly running mean is shown. Index computed by Enfield \emph{et al.} (2001) \cite{Enfield2001} using the Kaplan SST dataset provided by the NOAA/OAR/ESRL PSD, Boulder, Colorado, USA, from their website at https://www.esrl.noaa.gov/psd/.}
\end{figure}

The AMO is a pattern of variability in the  North Atlantic sea surface temperature  with a dominant period of fifty to seventy years \cite{Enfield2001}. In Figure \ref{fig:amoosc},  an index for  the average  sea surface temperature deviations in the North Atlantic Ocean is shown over the last 160 years. Although the instrumental record is somewhat limited for identifying the dominant time scale and  spatial pattern \cite{Deser2010}, such variability has been detected in proxy data \cite{Chylek2011} and in global climate models \cite{Delworth2000, Frankcombe2010b, Han2016, Cheung2017}. In most theories of the AMO \cite{DijkstraBook}, variations in the Atlantic Ocean circulation play a major role. The relevant  component of this circulation is the Atlantic Meridional Overturning Circulation (AMOC) \cite{Srokosz2012,DijkstraBookOcean}, which is basically the zonally averaged volume transport. The Gulf Stream is part of the AMOC, transporting warm water northward (and eastward) which loses heat on its way. At high northern latitudes the relatively heavy water sinks and  flows southward at larger depths.

One of the proposed physical mechanisms for the AMO, as described in Section \ref{sec:amo} below, is based on the propagation of so-called thermal Rossby waves \cite{Colin1999, TeRaa2002}. The role of these waves in the AMO motivates us to investigate whether its dynamics can be described by a delay model, since the propagation of waves underlies the delay in the ENSO model studied in \cite{Falkena2019}. Here, we thus apply the MZ formalism to a PDE model of the AMO, using a procedure inspired by \cite{Falkena2019}, with the aim of deriving a delay model describing the same AMO dynamics. The MZ formalism is general, can be applied to all types of equations (including non-linear and non-hyperbolic systems), and provides a more formal justification for the use of a conceptual delay model in analysing climate phenomena. This means the method discussed can be applied to other models in which there is a physical mechanism that could cause delayed effects.

The aim of the following study is two-fold: to derive a conceptual model with delay for the AMO and to demonstrate the utility and accuracy of doing this with the MZ formalism. We start with a brief description of the MZ formalism in Section \ref{sec:mz}, tailored to the problem at hand. Section \ref{sec:amo} presents the PDE model of the AMO by S{\'e}vellec and Huck \cite{Sevellec2015} which is the starting point for our study. The application of the MZ formalism to this model is discussed in Section \ref{sec:mzamo} where the resulting model with delay is introduced, as well as a model with delay that can be obtained for the simple PDE system using the method of characteristics. This leads to a comparison of the delay models and a discussion on the errors introduced by using the MZ formalism in Section \ref{sec:amodel}. A summary and discussion follows in Section \ref{sec:disc}.

\section{Mori-Zwanzig Formalism}
\label{sec:mz}

The Mori-Zwanzig (MZ) formalism provides a way of reducing a high-dimensional model to a reduced-order, more tractable system. The formalism is based on the work by Mori \cite{Mori1965} and Zwanzig \cite{Zwanzig1973} in statistical mechanics. It has been reformulated to be suitable for constructing reduced-order models for systems of ODEs \cite{Chorin2000, Chorin2002, Darve2009}. When studying applications in climate one often has to deal with PDEs, making the application of the formalism challenging. In this section we start with a general overview of the MZ formalism based on \cite{Chorin2002}, followed by a more detailed discussion of the particulars of applying the formalism to PDEs.

For a vector of state variables $\phi(x,t) \in \mathbb{R}^n$ which are continuously differentiable in $t \in \mathbb{R}_+$ and initial conditions $x \in \mathbb{R}^n$, we consider the system of ODEs defining the dynamics:
\begin{equation}
\label{eq:MZ_ODE}
\frac{\d}{\d t} \phi(x,t) = R( \phi(x,t) ), \qquad \phi(x,0) = x,
\end{equation}
where $R: \mathbb{R}^n \rightarrow \mathbb{R}^n$ is the vector-valued function of the specific system with components $R_i$. Now consider the evolution of an observable $u(x,t):=g(\phi(x,t))$ along a trajectory $\phi:\mathbb{R}^n \times \mathbb{R}_+ \rightarrow \mathbb{R}^n$. This observable satisfies the PDE
\begin{equation}
\label{eq:MZ_liouville_eq}
\frac{\partial}{\partial t} u(x,t) = \mathcal{L} u(x,t), \qquad u(x,0) = g(x),
\end{equation}
where $\mathcal{L}$ is the Liouville operator (or generator) \cite{morriss2013} given by
\begin{equation}
\label{eq:MZ_liouville_op}
[\mathcal{L}u](x) = \sum_{i=1}^n R_i(x) \partial_{x_i} u(x).
\end{equation}
Note that for a linear system where $R$ is defined by a matrix $A$ having elements $A_{ij}$, the Liouville operator reads $[\mathcal{L}u](x) = \sum_{i=1}^n \sum_{j=1}^n A_{ij} x_j \partial_{x_i} u(x)$.

To arrive at a reduced-order model for the dynamics governing $\phi(x,t)$, one needs to decide on the resolved variables $\hat{\phi}\in\mathbb{R}^m$. In our illustration we take $\hat{\phi}$ as a subset of components $\phi_i$ for some indices $i$. We also make a choice for an appropriate projection operator $P:C(\mathbb{R}^n, \mathbb{R}^k) \rightarrow C(\mathbb{R}^m, \mathbb{R}^k)$ onto these variables.
Examples of projection operators are the linear projection, setting all unresolved variables to zero, and the conditional expectation \cite{Chorin2002}. Let $Q = I - P$ denote the complement of $P$ (with $I$ the identity operator). Furthermore, we use the notation $[PR_i](\phi(x,t)) = R_i([\hat{\phi}(x,t),0]) = R_i(\hat{\phi}(x,t))$. We consider the choice of setting unresolved variables to zero for our projection $P$, such that, for an arbitrary observable $g\in C(\mathbb{R}^n, \mathbb{R}^k)$ (with arbitrary $k\geq1$), the projection $P$ is defined as $[Pg](\phi(x,t)) := g([\hat{\phi}(x,t),0])$.

Having chosen a set of resolved variables and a projection operator $P$, the reduced-order model corresponding to the full system \eqref{eq:MZ_ODE} is given by the generalized Langevin equation (see Chorin \emph{et al.} for its derivation \cite{Chorin2002}):
\begin{equation}
\label{eq:MZ_langevin}
\frac{\partial}{\partial t} \phi_i(x,t) = R_i([\hat{\phi}(x,t),0]) + F_i(x,t) + \int_0^t K_i([\hat{\phi}(x,t-s),0],s) \d s,
\end{equation}
where $\phi_i(x,t)$ is one of the resolved variables. The functions $F_i$ and $K_i$ are defined as
\begin{equation}
\label{eq:MZ_noise_memory}
F_i(x,t) = [e^{tQ\mathcal{L}}Q\mathcal{L}](x), \qquad K_i(\hat{x},t) = [P\mathcal{L}F_i](\hat{x},t),
\end{equation}
where $\hat{x}$ denotes the resolved part of the initial conditions $x$. The terms on the right-hand side of the Langevin equation are often referred to as the Markovian term $R_i(\hat{\phi}(x,t))$, the noise term $F_i(x,t)$ and the memory term, being the integral over $K_i(\hat{\phi}(x,t-s),s)$. For a linear system this memory integrand is obtained by applying a memory kernel to the resolved variables, i.e. $K_i(\hat{\phi}(x,t-s),s) = \hat{K}_i(s)[\hat{\phi}(x,t-s)]$.

The main difficulty in the application of the MZ formalism is calculating the terms $F_i$, which enter in the noise and the memory terms, and  which are the solutions of the orthogonal dynamics equations:
\begin{equation}
\label{eq:MZ_orthogonaldyn}
\frac{\partial}{\partial t} F_i(x,t) = Q\mathcal{L} F_i(x,t), \qquad F_i(x,0) = Q\mathcal{L} x_i.
\end{equation}
In general it is not known if the system \eqref{eq:MZ_orthogonaldyn} is well-posed. However, for specific cases it is possible to find approximate solutions. The possibility and difficulty of finding these solutions strongly depends on the choice of the resolved variables and projection operator. A suitable choice would yield an orthogonal dynamics system which can be solved in a more straightforward manner than the full system. In some cases,  the choice for the resolved variables and corresponding projection can be motivated by physical arguments for the specific system. For other models the choice might not be as straightforward and one cannot be certain that a suitable reduced-order model exists.

When the ODE system studied is linear, i.e. $R(\phi(x,t))=A\phi(x,t)$ where $A$ is a constant matrix, finding a suitable set of resolved variables can be done by looking at the eigenvalues of the orthogonal dynamics system. When the system \eqref{eq:MZ_ODE} is linear, the behaviour of the orthogonal dynamics can be obtained by studying the eigenvalues of $A_Q = A-PA$. A set of resolved variables is suitable if the eigenvalues of the orthogonal dynamics system show significantly more stability, i.e. have more negative real parts, compared to the full system. If this is not the case the problem of solving the full system is transferred to the equally difficult problem of solving the orthogonal dynamics system.

Up until now we have not discussed the difficulties arising when the MZ formalism is applied to a PDE system instead of ODEs. When the system is Hamiltonian some results exist (e.g. \cite{Chorin2000}), however, when this is not the case often the system is expanded in a basis of typically orthonormal functions (e.g. \cite{Darve2009, Zhu2018b}) to numerically find a solution. If the aim of applying the MZ formalism is to obtain a set of reduced-order model equations and not only a numerical result, this method is not suitable. Another approach, relying on integration along characteristics, has been explored by Falkena \emph{et al.} (2019) and yielded an exact reduced- order (delay) model for the system studied \cite{Falkena2019}. Here we build on this work to see whether delay type models can be derived for other systems of wave equations. In particular, we focus on a PDE model that describes thermal Rossby wave propagation related to the AMO.  This model is introduced in the next section.

\section{Atlantic Multidecadal Oscillation}
\label{sec:amo}

\begin{figure}
	\centering
	\begin{subfigure}{.45\textwidth}
		\centering
		\includegraphics[width=1.\textwidth]{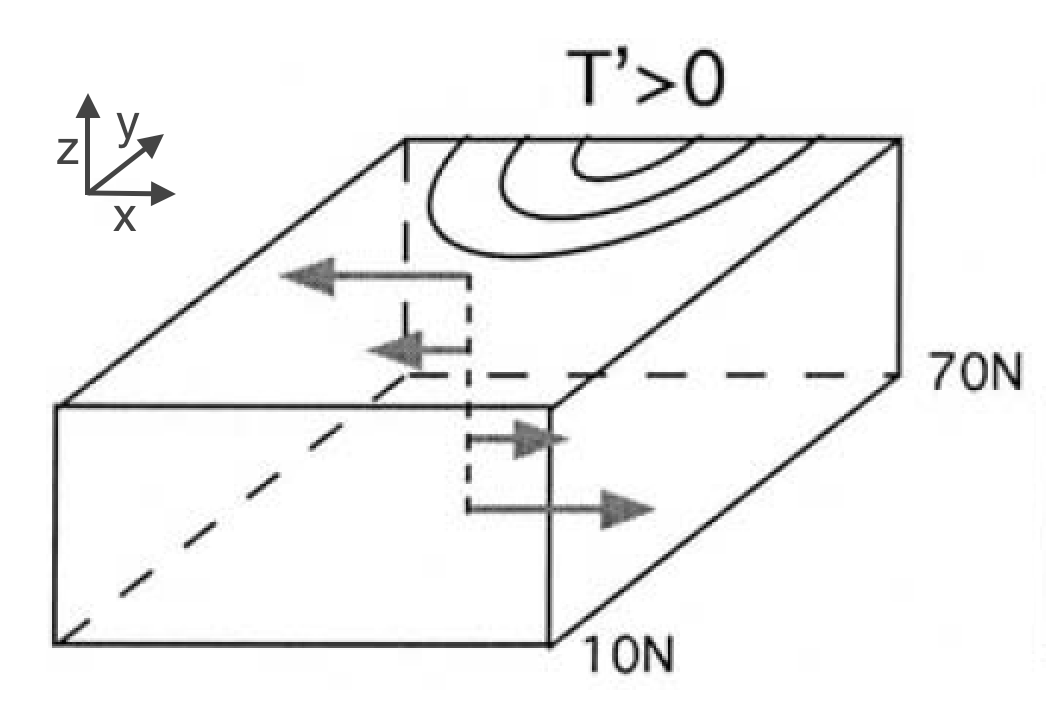}
		\caption{\label{fig:amophya}}
	\end{subfigure}
	\begin{subfigure}{.45\textwidth}
		\centering
		\includegraphics[width=1.\textwidth]{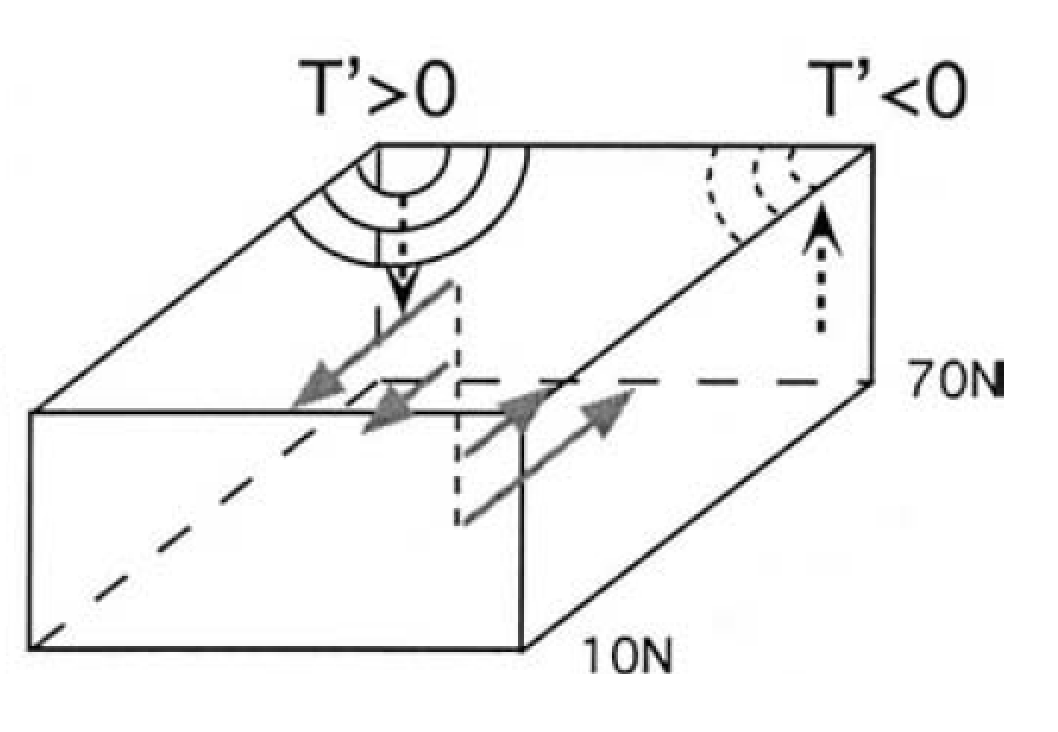} 
		\caption{\label{fig:amophyb}}
	\end{subfigure}
	\vspace{-10pt}
	\caption{\label{fig:amophy} Schematic diagram  of the physical mechanism responsible for the AMO with two phases a quarter period apart in panels (a) and (b). Figure taken from  \cite{Dijkstra2006}.
	}
\end{figure}

The thermal Rossby wave mechanism, suggested to be  responsible for the AMO \cite{Colin1999, TeRaa2002}, is  summarized in Figure \ref{fig:amophy}. When there is a positive temperature anomaly ($T'$) in the northern-central part of the basin, the meridional temperature gradient becomes stronger with respect to the background state. This results in a zonal overturning anomaly with westward surface flow through  thermal wind balance (Figure \ref{fig:amophya}). The negative zonal flow transports the positive temperature anomaly towards the western boundary, creating a zonal temperature gradient. Again through thermal wind balance, this now leads to anomalies in the meridional overturning circulation (Figure \ref{fig:amophyb}). This flow transports cold water from near the poles southward, reducing the meridional temperature gradient. This smaller north-south temperature gradient causes a positive (eastward) zonal flow, after which the same pattern as described above is followed with a sign change. Hence the variability associated with the AMO relies on the transport of heat and the flow response  through the thermal wind balance \cite{Dijkstra2006}.

Low-order ad-hoc ODE models of the AMO have been studied by, among others, Broer et al. \cite{Broer2011}. More recently, 
S\'evellec and Huck \cite{Sevellec2015} developed an idealized PDE model of the AMO model to which we apply the MZ 
formalism. This PDE model roughly captures the thermal Rossby wave mechanism described before, but with a 
simplification of the associated wave dynamics. In this section,  we present this model and briefly discuss its derivation. In 
addition we investigate the effect of using a more realistic strictly positive meridional overturning circulation as the background 
state on the behaviour of the \cite{Sevellec2015}  model.

\subsection{Model Formulation}
\label{ssec:formamo}

The AMO model by S\'evellec and Huck (2015) is a three-layer model describing the evolution of temperature perturbations in the North Atlantic Ocean \cite{Sevellec2015}. The model describes the temperatures ($T_i$, $i=1,2,3$) as a function of longitude ($x$) and time ($t$). For convenience we consider the non-dimensional version of the model with longitude-scale $W$ (basin width) and time-scale $Y$ (a year).

\begin{figure}
	\centering
	\includegraphics[width=.6\textwidth]{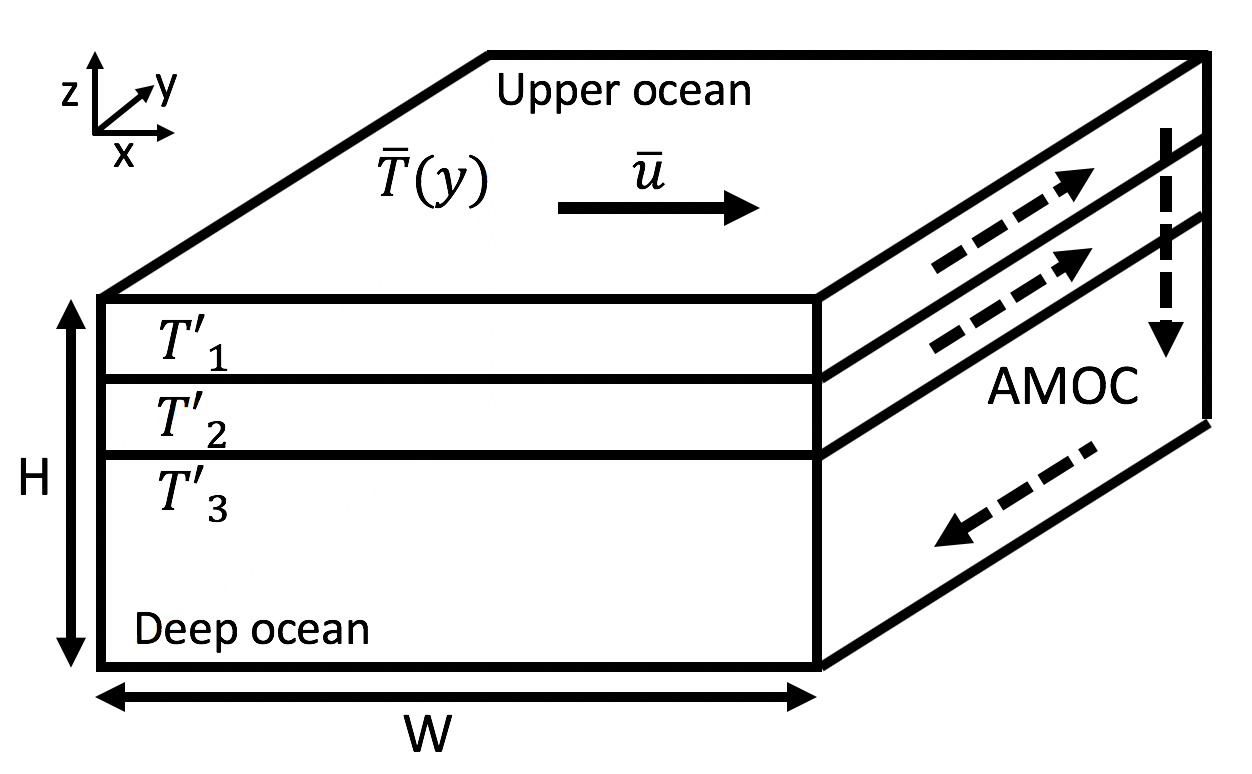}
	\vspace{-6pt}
	\caption{\label{fig:3layer_amoc} A schematic diagram of the three-layer ocean basin considered in the AMO-model. The dashed arrows show the background AMOC circulation, which is taken into account in the background state of the model discussed in Section \ref{sec:amo}\ref{ssec:backover}.}
\end{figure}

The scaled model is
\begin{equation}
\label{eq:modsh3}
\begin{split}
\partial_t T_1 & = a_1 \partial_x T_1 + b_1 \partial_x T_2 + c_1 \partial_x T_3 + \kappa_s \partial_{xx} T_1, \\
\partial_t T_2 & = a_2 \partial_x T_1 + b_2 \partial_x T_2 + c_2 \partial_x T_3 + \kappa_s \partial_{xx} T_2, \\
\partial_t T_3 & = \kappa_s \partial_{xx} T_3,
\end{split}
\end{equation}
with boundary conditions 
\begin{equation}
\label{eq:bcsh}
T_i(x|_{West}=0) = -T_i(x|_{East}=1), \qquad i=1,2,3.
\end{equation}
The constants in the model are all positive for physically realistic values and defined by
\begin{equation}
\label{eq:parval}
\begin{split}
a_1 &= \frac{Y}{W} \Big( \frac{\alpha_T g}{2Hf}\Big( -h_1(h_2+h_3) \partial_y \bar{T} + \frac{\beta}{2f} h_1^2(h_2+h_3)\partial_z \bar{T} \Big) - \bar{u}\Big), \\
b_1 &= \frac{Y}{W} \frac{\alpha_T g}{2Hf}\Big( -h_2(h_2+2h_3) \partial_y \bar{T} + \frac{\beta}{2f} h_1 h_2(h_2+2h_3)\partial_z \bar{T} \Big), \\
c_1 &= \frac{Y}{W} \frac{\alpha_T g}{2Hf}\Big( -h_3^2 \partial_y \bar{T} + \frac{\beta}{2f} h_1 h_3^2 \partial_z \bar{T} \Big),\\
a_2 &= \frac{Y}{W} \frac{\alpha_T g}{2Hf}\Big( h_1^2\partial_y \bar{T} + \frac{\beta}{2f} h_1^2(h_2+2h_3)\partial_z \bar{T} \Big), \\
b_2 &= \frac{Y}{W} \Big( \frac{\alpha_T g}{2Hf}\Big( -h_2(h_3-h_1)\partial_y \bar{T} + \frac{\beta}{2f} (4h_1h_2h_3 + h_2^2 (h_1+h_3))\partial_z \bar{T} \Big) - \bar{u} \Big), \\
c_2 &= \frac{Y}{W} \frac{\alpha_T g}{2Hf}\Big( -h_3^2\partial_y \bar{T} + \frac{\beta}{2f} h_3^2(2h_1+h_2)\partial_z \bar{T} \Big). \\
\kappa_s &= \kappa \frac{Y}{W^2}
\end{split}
\end{equation}
The values of the parameters needed to compute these constants are given in Table \ref{tab:parval}. 

\begin{table}
	\centering
	\caption{\label{tab:parval}The values of the parameters in the AMO model by S\'evellec and Huck \cite{Sevellec2015}.}
	\begin{tabular}{lll|c}
		\hline
		Thickness layer one & $h_1$ & 600 m & \\
		Thickness layer two & $h_2$ & 600 m & Vertical temperature gradient \\
		Thickness layer three & $h_3$ & 3300 m &  \\
		Total ocean depth & $H$ & 4500 m & $\partial_z \bar{T} = -\frac{2C}{h_1+h_2}\big(\Delta T - \frac{\alpha_S}{\alpha_T}\Delta S\big)$ \\
		Zonal basin size & $W$ & 4000 km &  \\
		Meridional basin size & $L$ & 6500 km & Control parameter $C$ \\
		Time scale (year) & $Y$ & $3.1536 \cdot 10^{7}$ s & Standard $C=1$ \\
		Horizontal diffusivity & $\kappa$ & $2 \cdot 10^3$ m$^2$/s &  \\
		Acceleration of gravity & $g$ & 9.8 m/s$^2$ &  \\
		Coriolis parameter & $f$ & $10^{-4}$ s$^{-1}$ \\
		$\beta$ effect & $\beta$ & $1.5\cdot 10^{-11}$ (ms)$^{-1}$ & Meridional temperature gradient \\
		Thermal expansion coefficient & $\alpha_T$ & $2\cdot 10^{-4}$ K$^{-1}$ & \\
		Haline contraction coefficient & $\alpha_S$ & $7\cdot 10^{-4}$ psu$^{-1}$ & $\partial_y \bar{T} = \frac{2}{L}\big(\Delta T - \frac{\alpha_S}{\alpha_T}\Delta S\big)$ \\
		Meridional temperature diff. & $\Delta T$ & -20 K &  \\
		Meridional salinity diff. & $\Delta S$ & -1.5 psu \\
		Zonal velocity & $\bar{u}$ & $10^{-2}$ m/s \\
		\hline
	\end{tabular}
\end{table}

The derivation of these equations can be found in \cite{Sevellec2015}. Here we briefly discuss that derivation and assumptions made to get to the above system of equations \eqref{eq:modsh3}. The derivation starts from an advection-diffusion equation for temperature, geostrophic balance (a balance between the horizontal pressure gradients and the Coriolis force), hydrostatic balance (a balance between the vertical pressure gradient and gravity) and a linear dependence of density on temperature. The equations for temperature are linearized around a fixed background state comprised of zonal flow $\bar{u}$ and temperature gradients in the meridional $\partial_y \bar{T}$ and vertical $\partial_z \bar{T}$. Note that this means that the overturning circulation ($\bar{v}$, $\bar{w}$) is neglected because of its weakness with respect to the zonal flow. The linearized temperature equation is then discretized over three layers assuming no flow through the surface and bottom and no background flow or temperature gradients in the bottom layer, which results in system \eqref{eq:modsh3}.

As the model assumes geostrophic balance, \eqref{eq:modsh3} only describes the solution of the interior part of the basin. The boundary conditions \eqref{eq:bcsh} are therefore derived by considering an additional boundary layer at either end of the basin with free-slip conditions at the interface between the interior flow and the boundary, and zero heat flux assumptions at the outer edges of the boundary layer (ocean basin walls). The full derivation can be found in the appendix of \cite{Sevellec2015}. Since the boundary conditions at hand will prove essential for our results, we explain the physics behind the coupling between the two boundaries. A signal, in the form of a Rossby wave, arriving at the western boundary of the basin, travels South, along the equator and back up North in the form of a Kelvin wave. Since the timescale of Kelvin waves is much shorter than that of the Rossby waves present in the model, this adjustment is assumed to be instantaneous. This leads to the coupling of the two boundaries, and allows for waves to keep propagating through the basin. For the specifics on the change of sign we refer the reader to the derivation of the boundary conditions in \cite{Sevellec2015}. Using for example a no-flux boundary is not valid here, as it would assume geostrophic balance in the boundary layer and neglect these Kelvin waves, making the model no longer dynamically accurate or suitable to study the AMO dynamics.

In this paper two additional simplifications to the model \eqref{eq:modsh3} are made. Firstly, we note that the only term acting in the third layer is diffusion and the two top layers do not couple into it. As a result any perturbation in that layer eventually damps out. For this reason, and to simplify the mathematical treatment of the system, perturbations in the bottom layer are neglected (i.e. $T_3=0$). Secondly, we approximate the diffusion terms by linear damping with damping coefficient $\alpha$. The system \eqref{eq:modsh3} then simplifies to a two-layer system:
\begin{equation}
\label{eq:modsh2}
\begin{split}
\partial_t T_1 & = a_1 \partial_x T_1 + b_1 \partial_x T_2 - \alpha T_1, \\
\partial_t T_2 & = a_2 \partial_x T_1 + b_2 \partial_x T_2 - \alpha T_2.
\end{split}
\end{equation}
This is the AMO model to which we apply the MZ formalism. Also note that this temperature model explains changes in the overturning circulation as well, via thermal wind balance, the continuity equation and Sverdrup balance, which is discussed in the Supplementary Information. The parameter values used for the numerical results in the remainder of the section and coming sections are given in Table \ref{tab:parnr}. 

\begin{table}
	\centering
	\caption{\label{tab:parnr}The numerical values of the parameters in Equation \eqref{eq:modsh2}. We note that $\alpha$ is a free parameter of $\mathcal{O}(10^{-3})$ and throughout this manuscript we will use $\alpha=0$ to explore the undamped solutions of Equation \eqref{eq:modsh2}.}
	\begin{tabular}{ccccc}
		\hline
		$a_1$ & $a_2$ & $b_1$ & $b_2$ \\ 
		0.1479 & 0.0540 & 0.4187 & 0.2423 \\ 
		\hline
	\end{tabular}
\end{table}

\begin{figure}[h]
	\centering
	\begin{subfigure}{.48\textwidth}
		\centering
		\includegraphics[width= \textwidth]{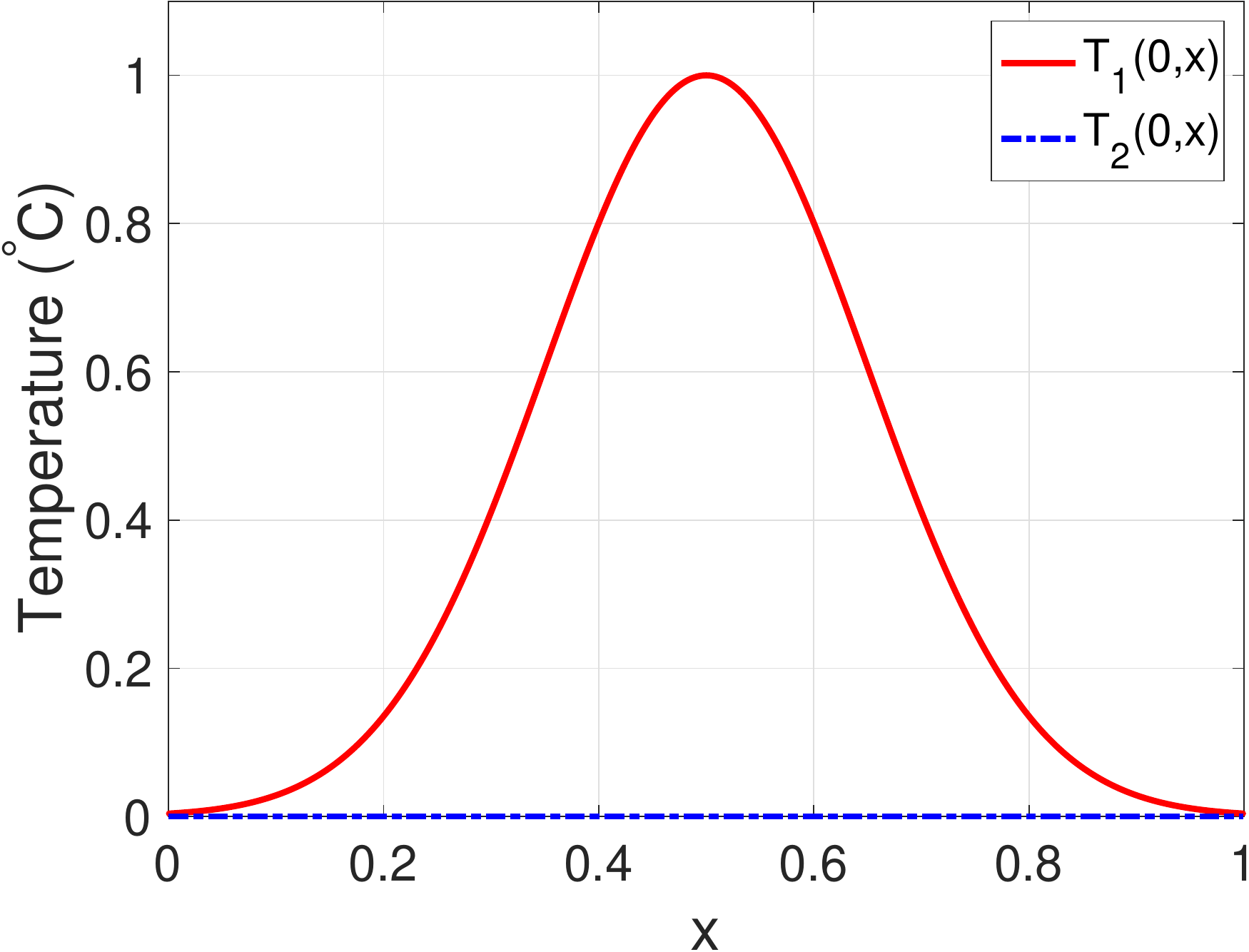}
		\caption{\label{fig:runsh_IC}The initial conditions.}
	\end{subfigure}
	\begin{subfigure}{.48\textwidth}
		\centering
		\includegraphics[width= \textwidth]{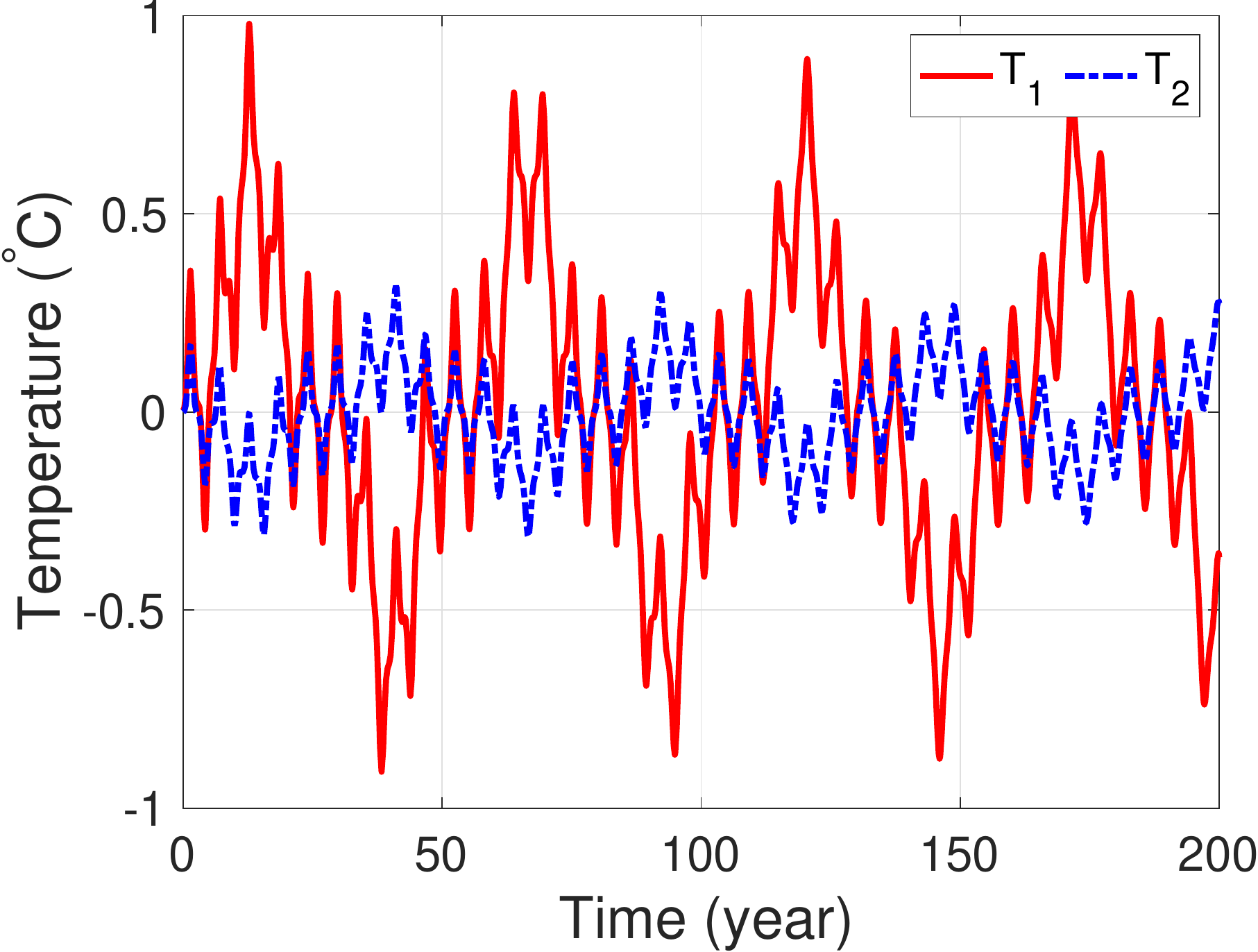}
		\caption{\label{fig:runsh_ev}Evolution of model \eqref{eq:modsh2} at $x=0$.}
	\end{subfigure}
	\caption{\label{fig:runsh}Model simulation of the temperature anomaly in the AMO model (Equation \eqref{eq:modsh2}) for an initial positive Gaussian temperature perturbation in the centre of the basin in the first layer ($\Delta t = \Delta x = 0.0005$, $\alpha=0$).}
\end{figure}

Before looking into the application of the MZ formalism to this AMO model, we illustrate its behaviour by simulating it for $\alpha=0$. We use an upwind discretization scheme for the $x$-derivatives and a forward Euler scheme in time. Note that this discretization includes numerical diffusion leading to artificial damping effects. The result is shown in Figure \ref{fig:runsh}. Note the opposite sign of the temperature in the two layers, which is due to the baroclinic nature of the waves \cite{GFDBook}. The model shows a combination of two oscillations with different periods. First, there is a long period of approximately sixty years, which corresponds to a thermal Rossby wave responsible for driving the AMO. Secondly, there is a higher frequency oscillation with a period of around five years.

The occurrence of the shorter period is at first sight surprising as it is not found in more detailed PDE models. This oscillation does not correspond to a planetary Rossby wave, as one might expect, since decreasing $\beta$ does not result in a disappearance of these oscillations.  It is a thermal Rossby wave, just as the one responsible for the AMO oscillation. The dominant appearance of this thermal Rossby wave in the model is  undesired to study the AMO. A possible improvement of the model, resulting in the damping of this high-frequency mode, is discussed in the next section. 

\subsection{Background Overturning Circulation}
\label{ssec:backover}

The AMO model by S\'evellec and Huck \cite{Sevellec2015} described in the previous section does not contain an overturning circulation in the background state, as the background meridional ($\bar{v}$) and vertical ($\bar{w}$) velocities are neglected. This means that in the model the overturning circulation, which can be inferred from the temperature evolution (details are given in the Supplementary Information), can become negative. To prevent this from happening in the model we consider an extended background state which retains meridional $\bar{v}$ and vertical $\bar{w}$ flow. With this different background state an extended two-layer temperature model for AMO can be derived following the same steps as in \cite{Sevellec2015}. The details of this derivation can be found in the Supplementary Information. The resulting system for temperature in the two upper layers is
\begin{equation}
\label{eq:modsh2_background}
\begin{split}
\partial_t T_1 & = a_1 \partial_x T_1 + b_1 \partial_x T_2 - (\beta_1 + \alpha) T_1 - \beta_2 T_2, \\
\partial_t T_2 & = a_2 \partial_x T_1 + b_2 \partial_x T_2 - (\beta_3 + \alpha) T_2,
\end{split}
\end{equation}
where
\begin{equation}
\label{eq:param_bg}
\beta_1 =  Y \cdot \Big( \frac{\beta}{f} \bar{v} + \frac{2}{h_1} \bar{w} \Big), \qquad \beta_2 = -Y \frac{4}{h_1} \bar{w}, \qquad \beta_3 =  Y \cdot \Big( \frac{\beta}{f} \bar{v} + \frac{2}{h_2} \bar{w} \Big).
\end{equation}
The difference to \eqref{eq:modsh2} is that there are additional linear terms in both equations. Note that not all the additional terms have a dampening effect, as some of the $\beta_i$-terms can be negative.

\begin{figure}[h]
	\centering
	\includegraphics[width=.5\textwidth]{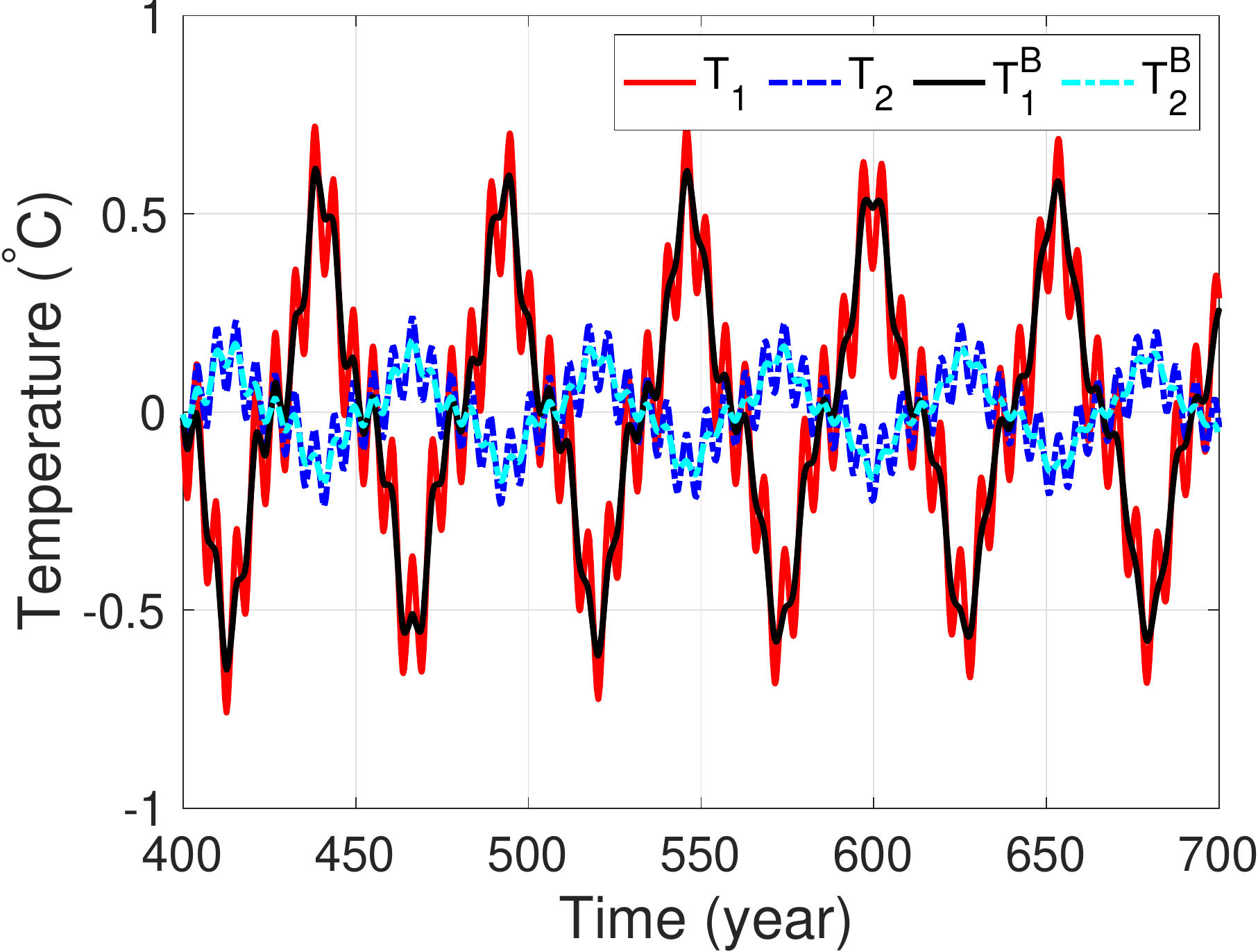}
	\caption{\label{fig:amobg}Simulations of the temperature in the two layers without (red, blue) and with (black, cyan) a background overturning circulation ($\beta_1 = \beta_3 = 1.156\cdot10^{-3}$, $\beta_2 = 7.148\cdot10^{-3}$, $\alpha=0$).}
\end{figure}

A model simulation for $\bar{v}=0.5\cdot 10^{-2} \text{ m/s}$ and $\bar{w}=-0.17\cdot 10^{-6} \text{ m/s}$ is shown in Figure \ref{fig:amobg}, where the values are chosen for plotting purposes being within a realistic range. Note that if $-2\bar{w}/h_{1,2} \gg \beta\bar{v}/f$, we have that $\beta_{1,3}$ become strongly negative leading to possible unstable solutions or at least amplifying effects within the solution. The result of adding the background overturning circulation is a damping of the high frequency oscillation, as can be seen in Figure \ref{fig:amobg}.
This short period oscillation is absorbed by the background overturning circulation while the long period oscillation persists. The amplitude of the oscillation corresponding to the AMO is not noticeably affected by the damping. This can be due to the presence of an amplifying effect of the background overturning in some parts of the equations as mentioned previously. For simplicity we apply the MZ formalism to the AMO model as given in Equation \eqref{eq:modsh2} instead of the extended model discussed here. The application of the MZ formalism to this extended model can be found in the Supplementary Material. This derivation follows the exact same steps as discussed in the next sections. We note that the results are similar to those discussed in the following, but for the extended model additional factors that lead to the decay of the high-frequency mode emerge. Thisis discussed in more detail at the end of Section \ref{sec:mzamo}\ref{ssec:amonm}.

\section{Reduction to a Delay Model}
\label{sec:mzamo}

The aim of applying the MZ formalism to the AMO model described in Section \ref{sec:amo} is to arrive at a projected model describing the same dynamics as the full model and analyze the effect of memory in that system, with the potential of deriving a delay model for the phenomenon.  A similar procedure has been applied to a model of the El Ni\~no Southern Oscillation (ENSO) by Falkena \emph{et al.} (2019) \cite{Falkena2019}. A difference is that for the AMO no previously proposed delay model is known. Therefore it is not immediately clear how to choose a projection, nor how to deal with solving the subsequent orthogonal dynamics equation \eqref{eq:MZ_orthogonaldyn}. Preferably we arrive at an equation for the temperature at one location in space, to remove the explicit dependency on $x$ in the system, but it is not clear from the onset whether or not this is feasible.

The way in which we proceed is to first convert the system of PDEs \eqref{eq:modsh2} into a set of ODEs by discretization. To this high-dimensional system of ODEs the MZ formalism is then applied. In the following sections this procedure is described. After deciding on the discretization to use, possible (sets of) resolved variables are explored, followed by a discussion of the different terms in the Langevin equation \eqref{eq:MZ_langevin}.

\subsection{Discretization}
\label{ssec:disc}

The first step is to find a stable discretization of the AMO model in Equation \eqref{eq:modsh2}. A grid of $(N+1)$-points in space with distance $dx = \frac{1}{N}$ is used. Because all parameters in the model are positive we know that all waves travel westward. Therefore we use an upwind scheme to discretize the model. The discretized equations are
\begin{equation}
\label{eq:amo_disc}
\begin{split}
\partial_t T_1^n &= \frac{a_1}{dx}(T_1^{n+1}-T_1^n) + \frac{b_1}{dx}(T_2^{n+1}-T_2^n) - \alpha T_1^n, \\
\partial_t T_2^n &= \frac{a_2}{dx}(T_1^{n+1}-T_1^n) + \frac{b_2}{dx}(T_2^{n+1}-T_2^n) - \alpha T_2^n,
\end{split}
\end{equation}
for $n=0,...,N$ (such that $T_i^k\approx T_i(k/N)$ and $dx=1/N$), with boundary conditions
\begin{equation}
T_1^{N} = -T_1^0,  \qquad T_2^{N} = -T_2^0.
\end{equation}
By the circular nature of the boundary conditions this is a $2N$-dimensional system (there is no need to solve the dynamical equations for discretization points $N$). Letting $N\rightarrow\infty$ recovers the PDE model exactly. This system \eqref{eq:amo_disc} can be written as a matrix equation for $\vec{T} = (T_1^0, T_2^0, ..., T_1^{N-1}, T_2^{N-1})$:
\begin{equation}
\partial_t \vec{T} = M \vec{T}.
\end{equation}
The construction of $M$ is straightforward from system \eqref{eq:amo_disc}.

\begin{figure}[h]
	\centering
	\begin{subfigure}{.48\textwidth}
		\centering
		\includegraphics[width= \textwidth]{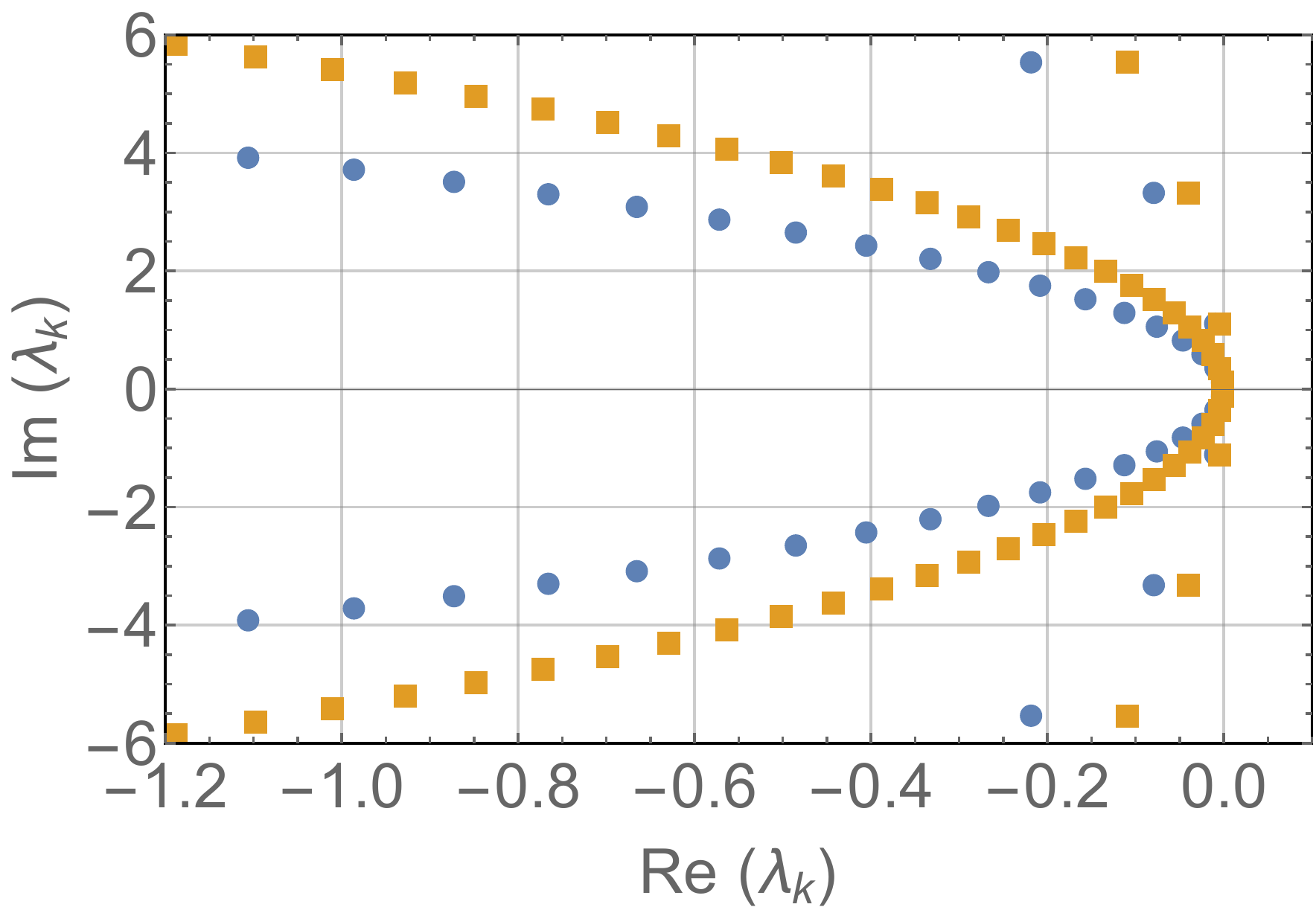}
		\caption{\label{fig:amoev_disc}Full model for $N=200$ (blue, circles) and $N=400$ (yellow, squares). \\}
	\end{subfigure}
	\hspace{10pt}
	\begin{subfigure}{.48\textwidth}
		\centering
		\includegraphics[width= \textwidth]{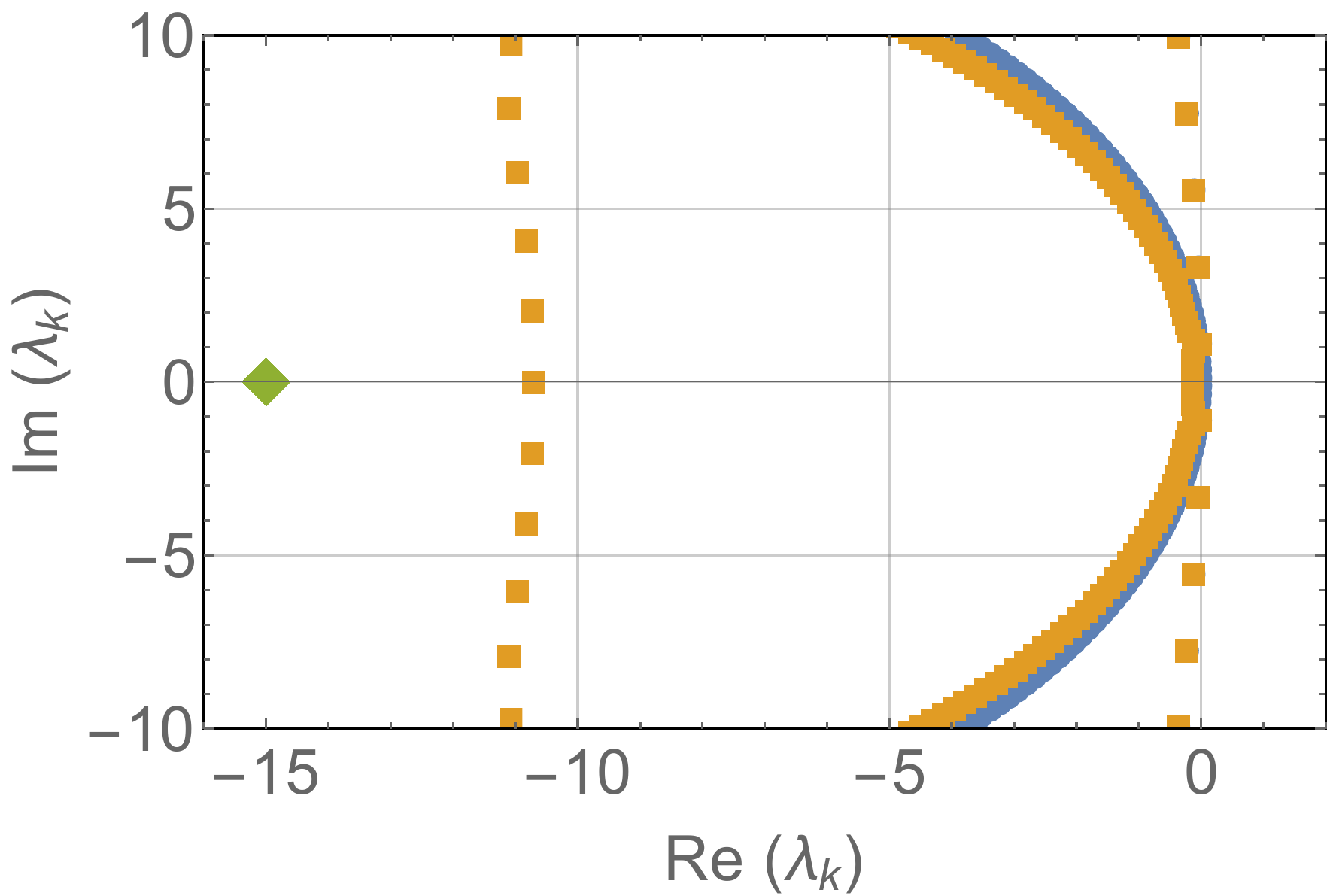}
		\caption{\label{fig:amoev_proj}Full model for $N=400$ (blue, circles) and using a projection onto either $T_1^0$ (yellow, squares) or both $T_1^0$ and $T_2^0$ (green, diamond).}
	\end{subfigure}
	\caption{\label{fig:amoev}The eigenvalues of the discretized AMO model ($\alpha=0.001$) (a) for different $N$ and (b) different sets of resolved variables. Note that not all eigenvalues are shown (e.g. only one of the two eigenvalues for the projection onto $T_1^0$ and $T_2^0$, the second one is more negative).
	}
\end{figure}

The stability of the solution of this discretized system of ODEs is verified by computing the eigenvalues of the matrix $M$. These are shown in Figure \ref{fig:amoev_disc} for $N=200$ (blue circles) and $N=400$ (yellow squares). For each $N$ two sets of eigenvalues are visible, with the spacing between the imaginary part of the eigenvalues in either set equal to the corresponding wave frequency. For increasing $N$ both curves of eigenvalues approach a line with real part $-\alpha$, being the eigenvalue of the continuous system. Since all eigenvalues are negative for every $N$ the discretization is stable. In the following sections we go into the application of the MZ formalism to this system of ODEs.

\subsection{Projection}
\label{ssec:mzproj}

When applying the MZ formalism the first step is to choose the resolved variables and corresponding projection. As discussed in Section \ref{sec:mz} this is an essential choice determining the final expressions in the Langevin equation \eqref{eq:MZ_langevin}. From the modelling perspective the aim is to find a system of equations for the temperature at one location in order to remove the $x$-dependency of the system. With this in mind there are three possible choices for the resolved variables; $T_1$, $T_2$, and both. Note that because all waves travel in the same direction without loss of energy it does not matter on which location in space the temperatures are projected. For convenience we choose to project onto the western boundary ($n=0$), but note that the result for any other location is the same. The most straightforward way to project onto one of those sets of resolved variables is to use the linear projection $P$, setting all unresolved variables to zero. Its complement $Q$ thus sets the resolved variables to zero since the system considered is linear. The corresponding orthogonal dynamics equation \eqref{eq:MZ_orthogonaldyn} can be written as a matrix equation
\begin{equation}
\label{eq:amo_ortdyn}
\partial_t \vec{T}_Q = M_Q \vec{T}_Q,
\end{equation}
where $\vec{T}_Q = (T_{1Q}^1, T_{2Q}^1, ..., T_{1Q}^{N-1}, T_{2Q}^{N-1})$ represents the unresolved variables. The matrix $M_Q$ is the same as the matrix for the full system $M$, but with the rows and columns corresponding to the resolved variable(s) removed.

To assess the quality of the three possible projections we study the eigenvalues of the corresponding orthogonal dynamics systems. A sufficient decay of the eigenvalues of the orthogonal dynamics system indicates the corresponding resolved variables and projections are suitable as discussed in Section \ref{sec:mz}. In Figure \ref{fig:amoev_proj} the eigenvalues of the full system and projected systems are shown for two different projections on boundary variables. When one projects onto only $T_1^0$ the eigenvalues are quite similar to those of the original system. Similar results are found for projection onto only $T_2^0$.

If both $T_1^0$ and $T_2^0$ are taken as the resolved variables the result is noticeably better. The orthogonal dynamics system has only two eigenvalues which become increasingly negative with increasing $N$. Therefore this projection onto both $T_1$ and $T_2$ at the boundary $x=0$ is chosen, where we note once again that the choice of the specific location is arbitrary. In the following section we focus on the derivation of the noise and memory term in the MZ formalism \eqref{eq:MZ_noise_memory}. Here we briefly discuss the Markovian terms. They are given by the projection of the right-hand side of the equation for $T_1^0$ and $T_2^0$ (Equation \eqref{eq:amo_disc} for $n=0$):
\begin{equation}
\label{eq:amo_markov}
\begin{split}
P\Big[\frac{a_1}{dx}(T_1^{1}-T_1^0) + \frac{b_1}{dx}(T_2^{1}-T_2^0) - \alpha T_1^0\Big] &= -\frac{a_1}{dx}T_1^0 - \frac{b_1}{dx}T_2^0 - \alpha T_1^0, \\
P\Big[\frac{a_2}{dx}(T_1^{1}-T_1^0) + \frac{b_2}{dx}(T_2^{1}-T_2^0) - \alpha T_2^0\Big] &= -\frac{a_2}{dx}T_1^0 - \frac{b_2}{dx}T_2^0 - \alpha T_2^0.
\end{split}
\end{equation}
This simply is the right-hand side dependence of Equation \eqref{eq:amo_disc} on the resolved variables, as a linear projection is used.

\subsection{Noise and Memory Terms}
\label{ssec:amonm}

To compute the noise and memory term we focus on the orthogonal dynamics system \eqref{eq:amo_ortdyn}. For the chosen resolved variables $T_1^0$ and $T_2^0$ with the linear projection the matrix reads
\begin{equation}
\label{eq:amo_MQ}
M_Q=
\begin{pmatrix}
-\frac{a_1}{dx}-\alpha & -\frac{b_1}{dx} & \frac{a_1}{dx}& \frac{b_1}{dx} & & \\
-\frac{a_2}{dx} & -\frac{b_2}{dx}-\alpha & \frac{a_2}{dx} & \frac{b_2}{dx} & & \\
& & -\frac{a_1}{dx}-\alpha & -\frac{b_1}{dx} & \ddots & & & \\
& & -\frac{a_2}{dx} & -\frac{b_2}{dx}-\alpha & & \ddots & & \\
& & & & \ddots & & \frac{a_1}{dx} & \frac{b_1}{dx} \\
& & & & & \ddots & \frac{a_2}{dx} & \frac{b_2}{dx} \\
& & & & & & -\frac{a_1}{dx}-\alpha & -\frac{b_1}{dx} \\
& & & & & & -\frac{a_2}{dx} & -\frac{b_2}{dx}-\alpha
\end{pmatrix}.
\end{equation}
Note this matrix is block upper diagonal with all blocks on the diagonal being the same. To solve the orthogonal dynamics system we have to find the eigenvalues and (generalized) eigenvectors of this matrix $M_Q$. As discussed in Section \ref{ssec:mzproj} there are only two eigenvalues
\begin{equation}
\label{eq:amo_eval}
\lambda_\pm = -\alpha - \frac{l_\pm}{dx},
\end{equation}
with
\begin{equation}
\label{eq:amomzl}
l_\pm = \frac{1}{2}\Big(a_1+b_2 \pm \sqrt{a_1^2+b_2^2-2a_1b_2+4a_2b_1}\Big),
\end{equation}
each of multiplicity $N-1$. Note that $l_\pm$ yield the characteristics of the original PDE system \eqref{eq:modsh2}. The corresponding generalized eigenvectors for $i = 1,...,N-1$ are
\begin{equation}
\label{eq:amo_evec}
\vec{v}_\pm^i = \Big( \frac{dx}{l_\pm} \Big)^{i-1} \cdot (0,...,0,w_\pm,1,0,...,0),
\end{equation}
where the non-zero values are located on the coordinates corresponding to location $i$. Here
\begin{equation}
w_\pm = \frac{1}{2a_2}\Big(a_1 - b_2 \pm \sqrt{a_1^2+b_2^2-2a_1b_2+4a_2b_1}\Big).
\end{equation}
Having computed the eigenvalues and eigenvectors we can write down the solutions $\vec{T}_Q$ of the orthogonal dynamics equation (e.g. \cite{Zill1997}). Here we only note that $\vec{T}_Q$ is a linear combination of the eigenvectors, meaning it is relatively straightforward to identify the solution at one location. The full expressions are given in the Supplementary Information, together with the use of initial conditions to determine the constants involved.

Now that we have the solution to the orthogonal dynamics equation we can write down the noise terms and subsequently compute the memory terms of the discretized AMO system \eqref{eq:amo_disc}. The noise terms are defined by
\begin{equation}
\begin{split}
F_{T_1^0}(t) &= \frac{a_1}{dx} T_{1Q}^1(t) + \frac{b_1}{dx} T_{2Q}^1(t), \\
F_{T_2^0}(t) &= \frac{a_2}{dx} T_{1Q}^1(t) + \frac{b_2}{dx} T_{2Q}^1(t).
\end{split}
\end{equation}
Note that only the terms of the solution $\vec{T}_Q$ which contain the eigenvectors $\vec{v}_\pm^1$ contribute to the noise term, as all other eigenvectors have zeros in the direction of $T_{1Q}^1$ and $T_{2Q}^1$. The resulting expressions, following the solution of Equation \eqref{eq:amo_ortdyn}, are
\begin{equation}
\label{eq:amo_noise}
\begin{split}
F_{T_1^0}(t) &= N  \sum_{i=1}^{N-1} \Big( (a_1 w_+ + b_1) e^{\lambda_+t} c_+^{i} + (a_1 w_- + b_1) e^{\lambda_-t} c_-^{i} \Big) \frac{t^{i-1}}{(i-1)!}, \\
F_{T_2^0}(t) &= N \sum_{i=1}^{N-1} \Big( (a_2 w_+ + b_2) e^{\lambda_+t} c_+^{i} + (a_2 w_- + b_2) e^{\lambda_-t} c_-^{i} \Big) \frac{t^{i-1}}{(i-1)!},
\end{split}
\end{equation}
where
\begin{equation}
\label{eq:amomzconstants}
\begin{split}
c_+^{i} &= \Big(\frac{l_+}{dx}\Big)^{i-1} \cdot \frac{T_1^{i}(0)-w_-T_2^{i}(0)}{w_+-w_-}, \\
c_-^{i} &= -\Big(\frac{l_-}{dx}\Big)^{i-1} \cdot \frac{T_1^{i}(0)-w_+T_2^{i}(0)}{w_+-w_-},
\end{split}
\end{equation}
depend on the initial conditions of the unresolved variables.

To compute the memory terms (as defined in \eqref{eq:MZ_noise_memory}) we first look at the effect of applying the operator $P\mathcal{L}$ to each of the initial conditions. This is sufficient for computation of the memory terms because the noise terms \eqref{eq:amo_noise} are linear in the initial conditions. We find
\begin{equation}
\begin{split}
& P \mathcal{L} \big(T_1^1(0), T_2^1(0), ..., T_1^i(0),T_2^i(0), ..., T_1^{N-1}(0),T_2^{N-1}(0)\big) \\
&= P \big(..., \frac{a_1}{dx}(T_1^{i+1}(0)-T_1^i(0)) + \frac{b_1}{dx}(T_2^{i+1}(0)-T_2^i(0)) - \alpha T_1^i(0) \\
& \qquad \qquad \frac{a_2}{dx}(T_1^{i+1}(0)-T_1^i(0)) + \frac{b_2}{dx}(T_2^{i+1}(0)-T_2^i(0)) - \alpha T_2^i(0) ,...\big), \\
&= \big(0, ...,0, -\frac{a_1}{dx}T_1^0(0) -\frac{b_1}{dx}T_2^0(0), -\frac{a_2}{dx}T_1^0(0) -\frac{b_2}{dx}T_2^0(0)\big).
\end{split}
\end{equation}
We see that only terms that initially depend on $T_{1}^{N-1}(0)$ and $T_{2}^{N-1}(0)$ are non-zero after application of $P\mathcal{L}$. Combining this result with the noise term \eqref{eq:amo_noise} and replacing $dx$ by $\frac{1}{N}$, the memory integrand \eqref{eq:MZ_noise_memory} becomes
\begin{equation}
\label{eq:amo_memory}
\begin{split}
K_{T_1^0}((T_1^0(0),T_2^0(0),t) &= 
N^2 \frac{t^{N-2}}{(N-2)!} e^{-\alpha t} \Big( (l_+ N)^{N-2} e^{-l_+ N t} \big( A_{1+} T_1^0(0) + B_{1+} T_2^0(0) \big) \\
&\qquad + (l_- N)^{N-2} e^{-l_- N t} \big( A_{1-} T_1^0(0) + B_{1-} T_2^0(0) \big) \Big), \\
K_{T_2^0}((T_1^0(0),T_2^0(0),t) 
&= N^2 \frac{t^{N-2}}{(N-2)!} e^{-\alpha t} \Big( (l_+ N)^{N-2} e^{-l_+ N t} \big( A_{2+} T_1^0(0) + B_{2+} T_2^0(0) \big) \\
&\qquad + (l_- N)^{N-2} e^{-l_- N t} \big( A_{2-} T_1^0(0) + B_{2-} T_2^0(0) \big) \Big),
\end{split}
\end{equation}
with
\begin{equation}
\label{eq:amoconstants}
\begin{split}
A_{1+} = \frac{(a_1 w_+ + b_1)(-a_1+w_-a_2)}{w_+ - w_-}, &\qquad A_{1-} = \frac{-(a_1 w_- + b_1)(-a_1+w_+a_2)}{w_+ - w_-}, \\
B_{1+} = \frac{(a_1 w_+ + b_1)(-b_1+w_-b_2)}{w_+ - w_-}, &\qquad B_{1-} = \frac{-(a_1 w_- + b_1)(-b_1+w_+b_2)}{w_+ - w_-}, \\
A_{2+} = \frac{(a_2 w_+ + b_2)(-a_1+w_-a_2)}{w_+ - w_-}, &\qquad A_{2-} = \frac{-(a_2 w_- + b_2)(-a_1+w_+a_2)}{w_+ - w_-}, \\
B_{2+} = \frac{(a_2 w_+ + b_2)(-b_1+w_-b_2)}{w_+ - w_-}, &\qquad B_{2-} = \frac{-(a_2 w_- + b_2)(-b_1+w_+b_2)}{w_+ - w_-}.
\end{split}
\end{equation}

Now all components of the Langevin equation \eqref{eq:MZ_langevin}, being the Markovian terms \eqref{eq:amo_markov}, the noise terms \eqref{eq:amo_noise} and the memory integrands \eqref{eq:amo_memory}, are known. Thus we can write down the result of applying the MZ formalism to the discretized AMO system \eqref{eq:amo_disc}:
\begin{equation}
\label{eq:mz_ndep}
\begin{split}
\partial_t T_1^0 &= -a_1 N T_1^0 - b_1 N T_2^0 - \alpha T_1^0 \\
&\quad +  N e^{-\alpha t} \sum_{i=1}^{N-1} \Big( (a_1 w_+ + b_1) e^{-l_+Nt} c_+^{i} + (a_1 w_- + b_1) e^{-l_-Nt} c_-^{i} \Big) \frac{t^{i-1}}{(i-1)!}\\
&\quad + \int_0^t N^2 \frac{(t-s)^{N-2}}{(N-2)!} e^{-\alpha (t-s)} \Big( (l_+ N)^{N-2} e^{-l_+ N (t-s)} \big( A_{1+} T_1^0(s) + B_{1+} T_2^0(s) \big) \\
&\qquad \qquad + (l_- N)^{N-2} e^{-l_- N (t-s)} \big( A_{1-} T_1^0(s) + B_{1-} T_2^0(s) \big) \Big) \d s, \\
\partial_t T_2^0 &= -a_2 N T_1^0 - b_2 N T_2^0 - \alpha T_2^0  \\
&\quad + N e^{-\alpha t} \sum_{i=1}^{N-1} \Big( (a_2 w_+ + b_2) e^{-l_+Nt} c_+^{i} + (a_2 w_- + b_2) e^{-l_-Nt} c_-^{i} \Big) \frac{t^{i-1}}{(i-1)!} \\
&\quad + \int_0^t N^2 \frac{(t-s)^{N-2}}{(N-2)!} e^{-\alpha (t-s)} \Big( (l_+ N)^{N-2} e^{-l_+ N (t-s)} \big( A_{2+} T_1^0(s) + B_{2+} T_2^0(s) \big) \\
&\qquad \qquad + (l_- N)^{N-2} e^{-l_- N (t-s)} \big( A_{2-} T_1^0(s) + B_{2-} T_2^0(s) \big) \Big) \d s.
\end{split}
\end{equation}

This system depends on the discretization, or more precisely, on the number of points $N$. Ideally we would like to find the equations for the continuous model. Therefore the limiting behaviour as $N\rightarrow\infty$, or $1/N = \epsilon \rightarrow 0$, of the different terms is studied in the Supplementary Material. The equation obtained after taking this limit can be written as
\begin{equation}
\label{eq:mzres}
\begin{split}
\epsilon \frac{\d T_1}{\d t} &= -a_1 T_1(t) - b_1 T_2(t) + A_{1+} \tau_+ e^{-\alpha \tau_+} T_1\big( t-\tau_+ \big) + B_{1+} \tau_+ e^{-\alpha \tau_+} T_2\big( t-\tau_+ \big) \\
&\quad + A_{1-} \tau_- e^{-\alpha \tau_-} T_1\big( t-\tau_- \big) + B_{1-} \tau_- e^{-\alpha \tau_-} T_2\big( t-\tau_- \big) + \epsilon f_{\epsilon1}(t) + \mathcal{O}(\epsilon^2), \\
\epsilon \frac{\d T_2}{\d t} &= -a_2 T_1(t) - b_2 T_2(t) + A_{2+} \tau_+ e^{-\alpha \tau_+} T_1\big( t-\tau_+ \big) + B_{2+} \tau_+ e^{-\alpha \tau_+} T_2\big( t-\tau_+ \big) \\
&\quad + A_{2-} \tau_- e^{-\alpha \tau_-} T_1\big( t-\tau_- \big) + B_{2-} \tau_- e^{-\alpha \tau_-} T_2\big( t-\tau_- \big) + \epsilon f_{\epsilon2}(t) + \mathcal{O}(\epsilon^2),
\end{split}
\end{equation}
where
\begin{equation}
\label{eq:mzerror}
\begin{split}
f_{\epsilon1}(t) &=  - \alpha T_1(t) + A_{1+} \tau_+ e^{-\alpha \tau_+}g_{\epsilon+}\big(T_1\big) + B_{1+} \tau_+ e^{-\alpha \tau_+} g_{\epsilon+}\big(T_2\big) \\
&\qquad + A_{1-} \tau_- e^{-\alpha \tau_-} g_{\epsilon-}\big(T_1\big)+ B_{1-} \tau_- e^{-\alpha \tau_-} g_{\epsilon-}\big(T_2\big), \\
f_{\epsilon2}(t) &= - \alpha T_2(t) + A_{2+} \tau_+ e^{-\alpha \tau_+} g_{\epsilon+}\big(T_1\big) + B_{2+} \tau_+ e^{-\alpha \tau_+} g_{\epsilon+}\big(T_2\big)\\
&\qquad + A_{2-} \tau_- e^{-\alpha \tau_-} g_{\epsilon-}\big(T_1\big) + B_{2-} \tau_- e^{-\alpha \tau_-} g_{\epsilon-}\big(T_2\big),
\end{split}
\end{equation}
for
\begin{equation}
g_{\epsilon\pm}\big(T\big) = \frac{\tau_\pm^2}{2} \Big( \big((l_\pm+\alpha)^2-\frac{7}{6}l_\pm^2\big) T\big( t-\tau_\pm \big) + 2(l_\pm+\alpha) T'\big( t-\tau_\pm \big) + T''\big( t-\tau_\pm \big) \Big),
\end{equation}
with $\tau_\pm=1/l_\pm$. Note that we have dropped the superscript 0 in the notation, as the found system is valid at every location throughout the basin through a simple coordinate transformation. In \eqref{eq:mzres} the first two terms (without delay) in the equations for $T_{1,2}$ are the Markovian terms, while the terms including a delay result from the memory term. In \eqref{eq:mzerror} the $\alpha$-term comes from the Markovian part, while the terms including $g_{\epsilon\pm}$ can be attributed to the memory-term. This is the final result of the application of the MZ formalism to the AMO model \eqref{eq:modsh2} as an expansion in terms of order $\epsilon$. Letting $\epsilon\rightarrow 0$ a set of delay difference equations is found, giving the exact reduced model of the AMO.

When applying the MZ formalism to the extended AMO model as derived in Section \ref{sec:amo}\ref{ssec:backover}, the leading order terms change slightly. The $e^{-\alpha \tau_\pm}$ terms change to $e^{-(\alpha+l^1_\pm) \tau_\pm}$ with $l^1_\pm$ the additional first order term of the eigenvalues of the orthogonal dynamics system for the extended AMO model. This leads to additional damping of the high-frequency modes as $l^1_+$ is positive, i.e. reducing the effect of the terms with a short delay time $\tau_+$. On the other hand $l^1_-$ is negative weakening the damping of the low-frequency mode (and making it weakly unstable for $\alpha=0$). This corresponds to the observed weakening of the high-frequency modes as discussed in Section \ref{sec:amo}\ref{ssec:backover}.

\subsection{Delay model derived via wave characteristics}
\label{ssec:delchar}

The Mori-Zwanzig formalism is in principle applicable also when the coefficients $a_j$ and $b_j$ are space dependent. However, in that case it may not be possible to derive explicit expressions for delays and coefficients in \eqref{eq:mzres}. For spatially constant coefficients $a_j$ and $b_j$ and equal damping $\alpha$ in all components we may also derive the leading orders of \eqref{eq:mzres} by  integration along wave characteristics. This approach is similar to that taken in \cite{Falkena2019} and we refer to it as the method of characteristics (MoC). The damped free-wave solutions of the two-layer system,
\begin{equation}
\label{eq:free-wave}
\begin{split}
\partial_t T_1 & = a_1 \partial_x T_1 + b_1 \partial_x T_2-\alpha T_1\mbox{,} \\
\partial_t T_2 & = a_2 \partial_x T_1 + b_2 \partial_x T_2-\alpha T_2\mbox{,}
\end{split}
\end{equation}
solve the decoupled equations
\begin{align}
\label{eq:sys_tilde:pm}
\partial_t\tilde{T}_\pm-l_\pm\partial_x \tilde{T}_\pm+\alpha \tilde{T}_\pm&=0\mbox{,}
\end{align}
along the wave characteristics, where the characteristic speeds and delays
\begin{align} \label{eq:eigs_char}
l_\pm &= \frac{1}{2}\left[a_1+b_2 \pm \sqrt{(a_1 + b_2)^2 - 4a_1b_2+4a_2b_1}\right]\mbox{,}&
\tau_\pm&=1/l_\pm\mbox{,}
\end{align}
are the same as for the Mori-Zwanzig formalism in \eqref{eq:amomzl}. The new variables $\tilde{T}_\pm$ are related back to the original variables through the transformation $T_P$,
\begin{align} \label{eq:PT}
\vec{T}=&T_P
\begin{bmatrix}
\tilde{T}_+\\\tilde{T}_-
\end{bmatrix}
\mbox{,}& \mbox{where\quad} T_P &=
\begin{bmatrix}
\frac{a_1-l_+}{a_2} & \frac{a_1-l_-}{a_2}\\ 1&1
\end{bmatrix}\mbox{,}&
T_P^{-1} &= \begin{bmatrix} \frac{-a_2}{l_+-l_-} & \frac{a_1-l_-}{l_+-l_-} \\ \frac{a2}{l_+-l_-} & \frac{l_+-a_1}{l_+-l_-} \end{bmatrix}\mbox{.}
\end{align}
The damped wave equations \eqref{eq:sys_tilde:pm} have the general solutions
\begin{equation}\label{eq:char:genform}
\tilde{T}_\pm(t,x)=\tilde{T}_{\pm}^0(t+x\tau_\pm )\e^{-\alpha x\tau_\pm}\mbox{,}
\end{equation}
where the arbitrary profiles $\tilde{T}_{\pm}^0$ are constrained by the boundary conditions $\tilde{T}_\pm(t,0)=-\tilde{T}_\pm(t,1)$. These boundary conditions enforce the delay-difference equations for $\tilde{T}_\pm^0$,
\begin{equation} \label{eq:withBCs}
\tilde{T}_\pm^0(t) = -\e^{-\alpha\tau_\pm}\tilde{T}_\pm^0(t-\tau_\pm)\mbox{,}
\end{equation}
after shifting $t$ by $\tau_\pm$ and multiplying both sides by $\e^{-\alpha\tau_\pm}$ in the boundary conditions. Transforming $\tilde{T}^0(t)$ back using the transformation $T_P^{-1}$ gives the coupled delay equations
\begin{equation} \label{eq:2D-diff-vec}
\vec{T}(t) = \e^{-\alpha\tau_+}C_1\vec{T}(t-\tau_+) + \e^{-\alpha\tau_-}C_2\vec{T}(t-\tau_-),
\end{equation}
(dropping the superscript $0$ of $\tilde{T}$ for ease of notation) with
\begin{align*} \label{eq:2D-coeff}
C_1 = \begin{bmatrix} 
\frac{l_--a_1}{l_+-l_-} & 
\frac{(l_+-a_1)(l_--a_1)}{a_2(l_+-l_-)} \\  
-\frac{a2}{l_+-l_-} &  
-\frac{l_+-a_1}{l_+-l_-} 
\end{bmatrix}, &\qquad
C_2 = \begin{bmatrix} 
- \frac{l_+-a_1}{l_+-l_-} & 
- \frac{(l_+-a_1)(l_--a_1)}{a_2(l_+-l_-)} \\ 
\frac{a_2}{l_+-l_-} & 
\frac{l_+-a_1}{l_+-l_-}  \end{bmatrix}.
\end{align*}
The above is valid at every location in the basin according to the general solution form \eqref{eq:char:genform}, as discussed in Section \ref{sec:mzamo}\ref{ssec:amonm}. The MZ derived system in Equation \eqref{eq:mzres} for $\epsilon=0$ can be rewritten to the above system.

System \eqref{eq:2D-diff-vec} is a delay-difference system. In order to explore solutions of delay-difference system, we convert \eqref{eq:2D-diff-vec} to a system of DDEs by regularising it with a small time derivative $\epsilon\mathrm{d} \vec{T}/\mathrm{d}t$:
\begin{equation} \label{eq:DDE-AMO}
\epsilon \frac{d\vec{T}(t)}{dt}= -\vec{T}(t) + \e^{-\alpha\tau_+}C_1\vec{T}(t-\tau_+) + \e^{-\alpha\tau_-}C_2\vec{T}(t-\tau_-),
\end{equation}
with $\epsilon\ll 1$.  The choice of $\epsilon$ is related to the discretization of the original PDE system \eqref{eq:free-wave} through $\epsilon = 1/N$, where $N$ is the number of discretization steps using an `upwind' scheme (discretizing in the direction of the wave). Table \ref{tab:param3} shows the approximate resulting wave speeds and delays when using the parameter values in Table \ref{tab:parnr} for our numerical solutions and spectral analysis, which are discussed in the next section. 

\begin{table}
	\begin{center}
		\begin{tabular}{ | c | c | c | c |} 
			\hline
			$l_{+}$ & $l_{-}$ & $\tau_{+}$ &  $\tau_{-}$\\ 
			0.3527 & 0.0375 & 2.83 & 26.65 \\ 
			\hline
		\end{tabular}
	\end{center}
	\caption[Parameters for AMO model with delay]{Parameters used in numerical computations of (\ref{eq:DDE-AMO}). Note that these are approximated from Equation \eqref{eq:eigs_char} using the values in Table \ref{tab:parnr}. \label{tab:param3}}
\end{table}

\section{Analysis of Delay Models}
\label{sec:amodel}

In this section we analyze the solutions of the delay difference models derived via the MZ formalism \eqref{eq:mzres} and the MoC \eqref{eq:DDE-AMO}. We start with a discussion of the asymptotic spectrum and a spectral analysis of the MoC model \eqref{eq:DDE-AMO}. Next, we discuss the error terms as computed using the MZ formalism, followed by a comparison of the MZ delay model, the MoC delay model and the numerical PDE solution from which the MZ model is derived.

\subsection{Asymptotic Spectrum of Delay Models}
\label{sec:spec:asy}

We observe that the delays occurring in the MoC system are of different magnitude: $\tau_+\ll\tau_-$, where $\tau_+$ is of order $O(1)$ in the time scale of \eqref{eq:DDE-AMO}. For hierarchical large delays \cite{lichtner2011spectrum} and \cite{ruschel2019spectrum} provide a simple approximation of the spectrum for \eqref{eq:DDE-AMO} which captures the range of possible curvatures of the curves along which the eigenvalues shown in Fig.~\ref{fig:amoev} align. Any eigenvalue $\lambda$ of \eqref{eq:DDE-AMO} satisfies
\begin{equation}\label{spec:dde-amo}
\det\left[-\epsilon\lambda I-I+C_1\e^{-(\alpha+\lambda)\tau_+}+C_2\e^{-(\alpha+\lambda)\tau_-}\right]=0
\end{equation}
with $\epsilon\ll1$, $\tau_+=O(1)$ and $\tau_+\ll\tau_-$. Hence the
term $\e^{-(\alpha+\lambda)\tau_-}$ is negligibly small unless
$\alpha\tau_-$ and $\re\lambda\tau_-$ are of order $1$ or less. If we
assume that this is the case, we may introduce
$\alpha_-=\alpha/\tau_-$ (which is of order $1$ or less) and look for
eigenvalues $\lambda$ of the form
$\lambda=\gamma/\tau_++\i\omega/\epsilon$ (called the
pseudo-continuous spectrum in \cite{lichtner2011spectrum}). Then
$(\gamma+\alpha_-)\tau_+/\tau_-$ and $\epsilon\gamma/\tau_-$ are
small. Dropping these terms and introducing the phases
$\phi_\pm=\omega\tau_\pm/\epsilon$ and
$z=\e^{-(\gamma+\alpha_-)+\i\phi_-}$ simplifies \eqref{spec:dde-amo}
to
\begin{equation}
\label{eq:spec:asy}
\det\left[-\i\omega I-I+C_1\e^{\i\phi_+}+C_2z\right]=0  
\end{equation}
Equation~\eqref{eq:spec:asy} is in our case a quadratic equation in
the complex number $z$, giving two roots, each depending on $\omega$
and $\phi_+$, which one may express as $z_\pm(\omega,\phi_+)$. From
this root pair one may derive the damping depending on the frequency,
$\gamma_\pm(\omega,\phi_+)=-\alpha_-+\log z_\pm(\omega,\phi_+)$, and, in the original scaling for eigenvalue $\lambda$
\begin{displaymath}
\re\lambda=-\alpha+\tau_+\log z_\pm(\epsilon\im\lambda,\phi_+)\mbox{.}
\end{displaymath}
This relation determines the curves along which the eigenvalues align
for positive small $\epsilon$ and $\tau_+\ll\tau_-$. The phase
$\phi_+$ is treated here as an independent parameter. It is very
sensitive with respect to small changes of $\tau_+$ (since
$\phi_+=\omega\tau_+/\epsilon$) such that the location of the
eigenvalue curves will vary strongly depending on $\tau_+$ or
$\epsilon$ within the range given by $\phi_+\in[0,2\pi]$. Ruschel and
Yanchuk's \cite{ruschel2019spectrum} analysis shows in general that
for hierarchically large delays the spectrum ``fills an area'' of the
complex plane under small parameter variations.

\subsection{Spectral Analysis of Trajectories}
\label{sec:traj}
\begin{figure}
	\centering
	\begin{subfigure}{0.49\textwidth}
		\centering
		\includegraphics[width=\textwidth]{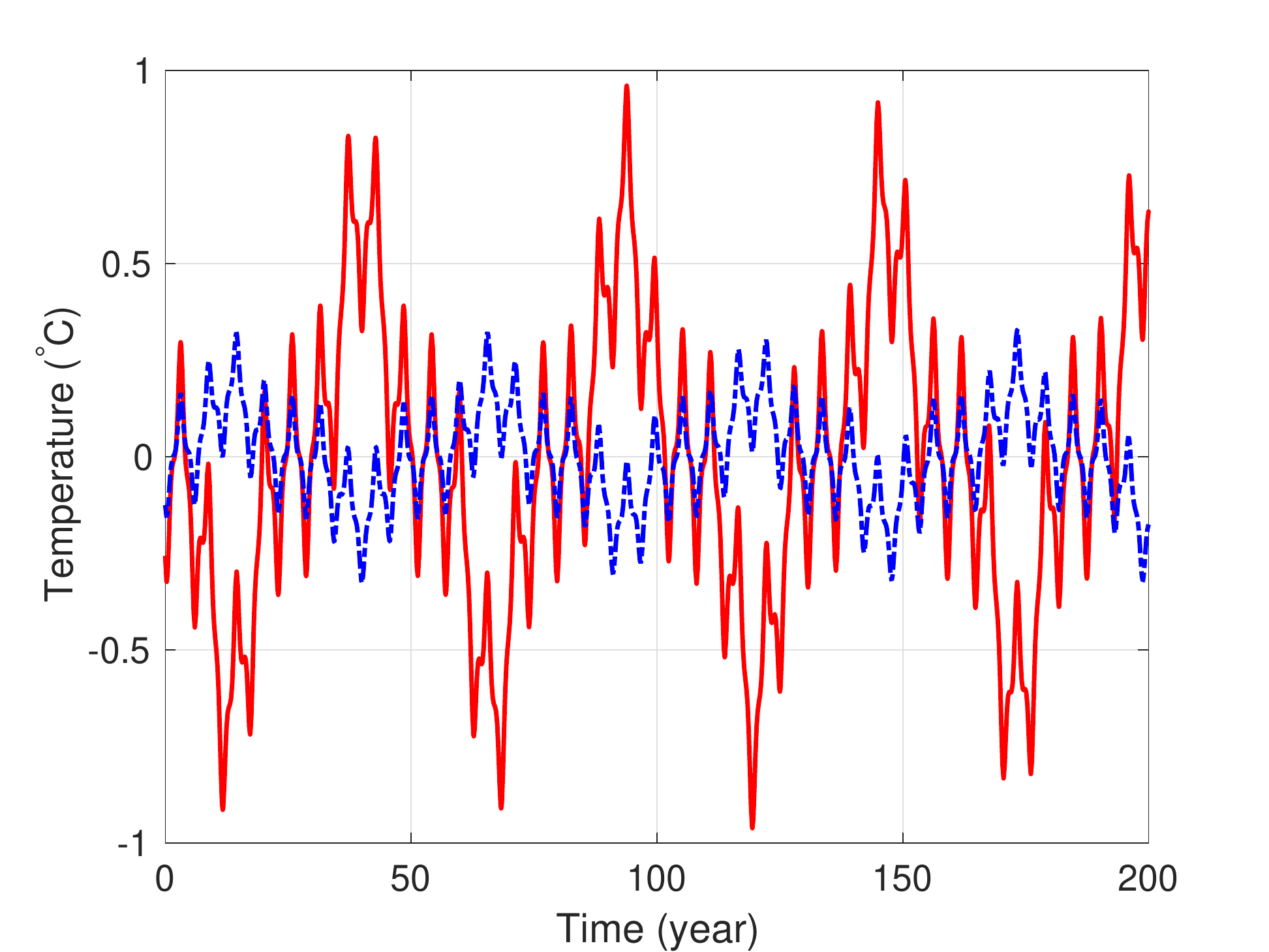}
		\caption{}
		\label{fig:traj_c08}
	\end{subfigure}
	\begin{subfigure}{0.49\textwidth}
		\centering
		\includegraphics[width=\textwidth]{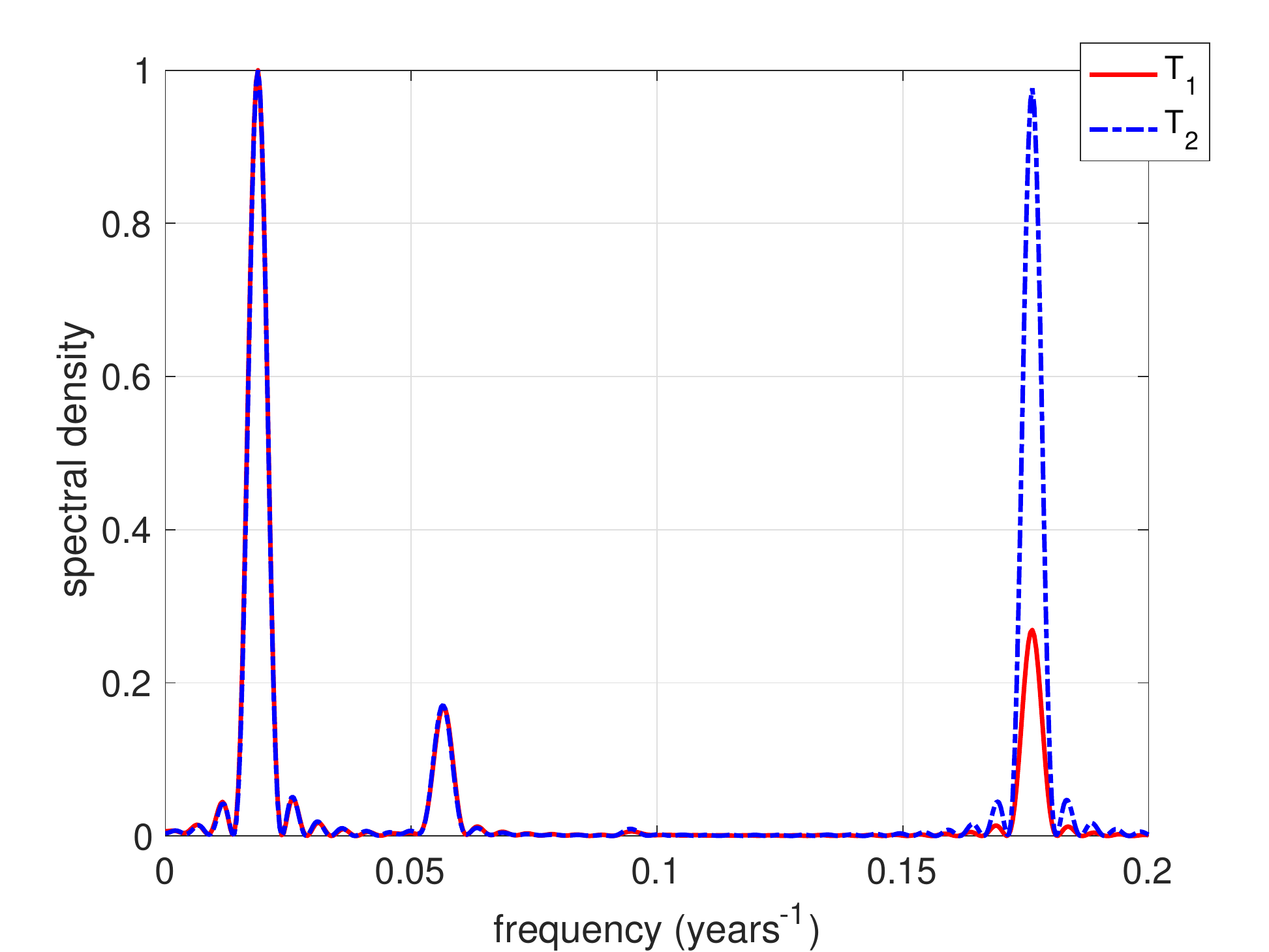}
		\caption{}
		\label{fig:spec_c08}
	\end{subfigure}
	\caption[Numerical results of DDE model for AMO]{Numerical results for (\ref{eq:DDE-AMO}) with parameters from Table \ref{tab:param3} and $\alpha=0$. (a)  Trajectory for 200 years. (b) Power spectral density.} 
	\label{fig:c08}
\end{figure}

For the spectral analysis of the MoC delay model \eqref{eq:DDE-AMO} (with $\alpha=0$) we compute the trajectories of $T_1$ and $T_2$. To compute the history needed for the difference equation, the PDE system (\ref{eq:free-wave}) is solved numerically for an initial profile of the basin using an upwind discretization scheme for $\tau_-$ years.  We take a Gaussian initial distribution profile (same as Figure \ref{fig:runsh_IC}).  The DDE system (\ref{eq:DDE-AMO}) is then evolved for a further 200 time steps.  Figure \ref{fig:c08} shows the results.

A spectral analysis is performed on the resulting trajectories to identify the most prominent oscillation periods.  A dominant signal of a $2\tau_-$ year cycle is obtained, along with a smaller signal for a $2\tau_+$ year cycle.  These two most prominent signals correspond to period doubling of the two delay values which arises naturally from the boundary conditions. There is a smaller peak corresponding to a cycle of approximately $\frac{2}{3}\tau_{-}$ years.  The signals corresponding to $2\tau_-$ (53.07) and $\frac{2}{3}\tau_{-}$ (17.77) year cycles align with the literature regarding possible cycle lengths of the AMO \cite{Delworth2000,Chylek2011}.  The signal corresponding to the $2\tau_+$ year cycle is much less pronounced in the surface temperature compared to the subsurface temperature.

The length of the cycles will naturally be dependent upon the basin size chosen.  Here we use a zonal basin size of $4000$km which corresponds to the width of the basin at $52^\circ$N.  The basin was defined in \cite{Sevellec2015} to have latitudinal boundaries of $10^\circ$N and $70^\circ$N.  The zonal basin width for a longitudinal extent of $60^\circ$ varies between these latitudes. We have therefore performed a sensitivity analysis to the basin width for a selection of latitudes and computed the resulting dominant cycle lengths.  The results of this are listed in Table \ref{tab:sens_basin}.  The range of cycle lengths, particularly those associated with $2\tau_-$ and $\frac{2}{3}\tau_-$, generally agree with the range of those identified in observational products. Measurements of SST have shown oscillations with a period between 50-70 years \cite{kushnir1994interdecadal,schlesinger1994oscillation}, while an analysis of sub-surface temperature identifies a 20-30 year oscillation \cite{frankcombe2008sub}. The first period of 50-70 years corresponds to $2\tau_-$, while the shorter period of 20-30 years roughly matches $\frac{2}{3}\tau_-$, with the best correspondence found for latitudes around 40$^\circ$N.

\begin{table}
	\centering
	\caption{\label{tab:sens_basin}Sensitivity analysis on basin length ($W$, km) and effect on dominant cycles (years) for $\alpha=0$.}
	\begin{tabular}{ccccc}
		\hline
		latitude & $W$& $2\tau_-$ & $\frac{2}{3}\tau_-$ &  $2\tau_+$ \\ \hline
		$10^\circ$N & $6540$ & $87.15$ & $29.05$ & $9.27$ \\
		$20^\circ$N & $6240$ & $83.15$ & $27.72$ & $8.85$ \\
		$30^\circ$N & $5760$ & $76.75$ & $25.58$ & $8.17$ \\
		$40^\circ$N & $5100$ & $67.96$ & $22.65$ & $7.23$ \\
		$50^\circ$N & $4260$ & $56.77$ & $18.92$ & $6.04$ \\
		$60^\circ$N & $3360$ & $44.77$ & $14.92$ & $4.76$ \\
		$70^\circ$N & $2280$ & $30.38$ & $10.13$ & $3.23$ \\
		\hline
	\end{tabular}
\end{table}

\subsection{Comparison of Delay Models}
\label{ssec:compdel}

\begin{figure}[h]
	\centering
	\begin{subfigure}{.48\textwidth}
		\centering
		\includegraphics[width= \textwidth]{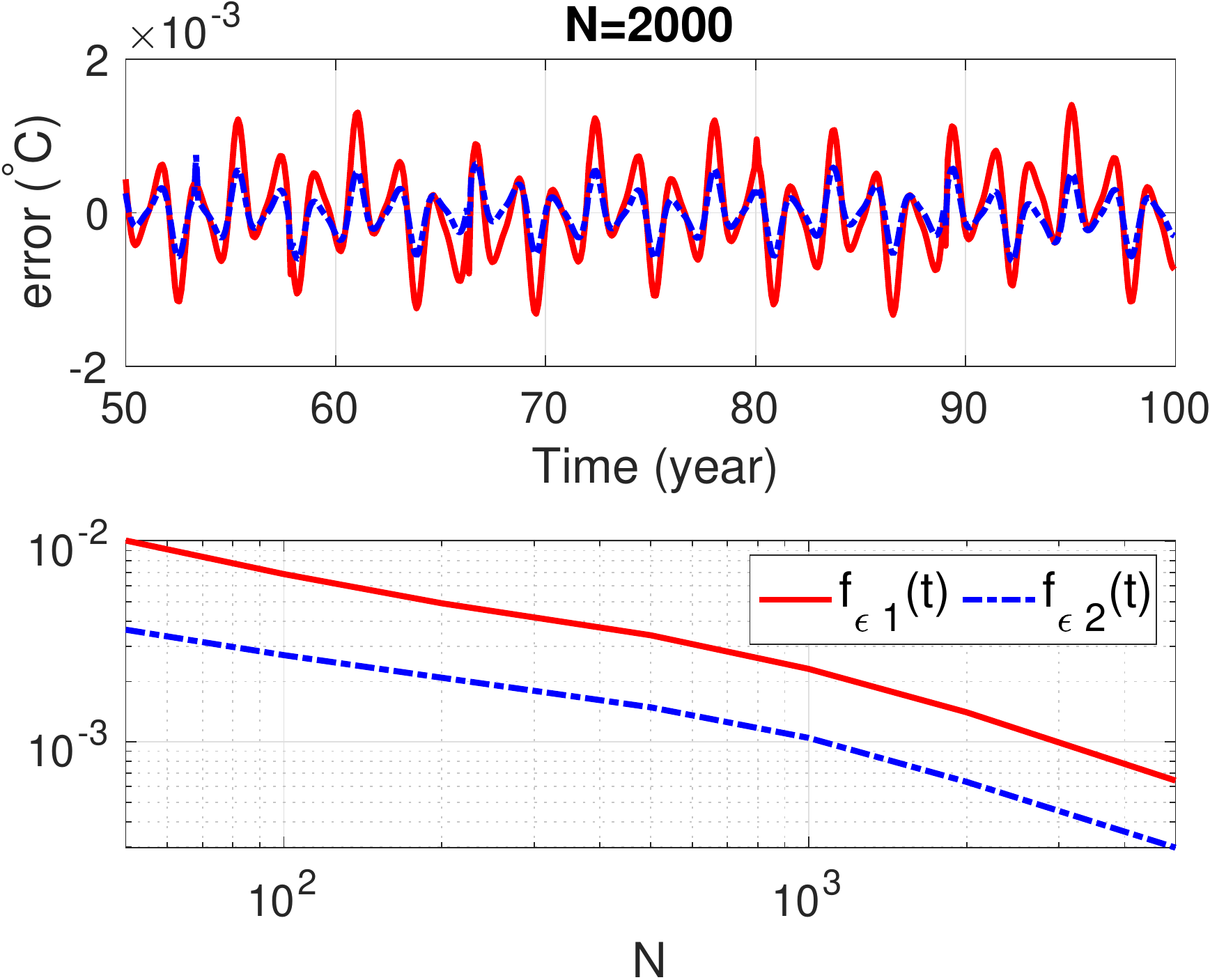}
		\caption{\label{fig:error_theo}}
	\end{subfigure}
	\begin{subfigure}{.48\textwidth}
		\centering
		\includegraphics[width= \textwidth]{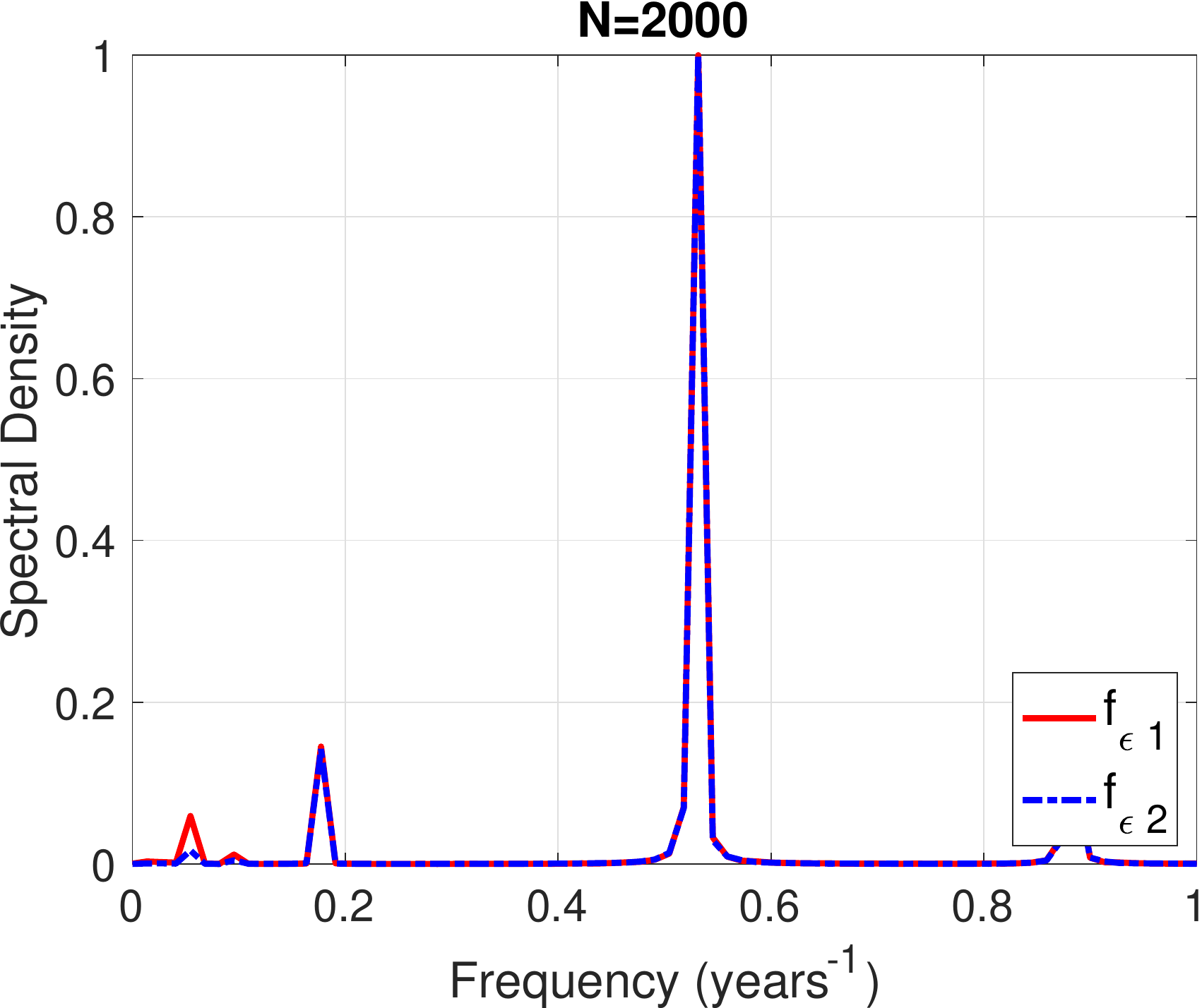}
		\caption{\label{fig:error_theo_freq}}
	\end{subfigure}
	\caption{\label{fig:deltheo}
		The error terms $f_{\epsilon 1}$, $f_{\epsilon 2}$ as computed using the MZ formalism with parameters from Table \ref{tab:param3} and $\alpha=0$. (a) The change of amplitude of the terms with $N$ is shown together with an example of the error terms for $N=2000$. (b) The spectrum of $f_{\epsilon 1}$ and $f_{\epsilon 2}$.}
\end{figure}

We start this section with an evaluation of the theoretical error terms $f_{\epsilon,i}$ for $i=1,2$ (Equation \eqref{eq:mzerror}) as computed via the MZ formalism. An example of the evolution of these terms has been plotted in Figure \ref{fig:error_theo}. Also shown is the decrease of their maximum amplitude with increasing $N$ (decreasing $\epsilon$), which corresponds to the effect expected from increasing  the number of steps in a discretized PDE. Thus the theoretical error of the delay system derived using the MZ formalism \eqref{eq:mzres} indicates that the delay model exhibits an error similar to that of the discretized PDE. In Figure \ref{fig:error_theo_freq} the spectral density of $f_{\epsilon i}$ for $i=1,2$ is shown. We find a peak at a frequency of ~0.53 years$^{-1}$ which emerges due to second derivatives in the computation of $f_{\epsilon i}$. The small low-frequency peaks correspond to those of the exact delay model, as shown in Figure \ref{fig:spec_c08}, as there is a delayed contribution from the temperature itself(and its first derivative) to the error term as well.

\begin{figure}
	\centering
	\begin{subfigure}{0.48\textwidth}
		\centering
		\includegraphics[width=\textwidth]{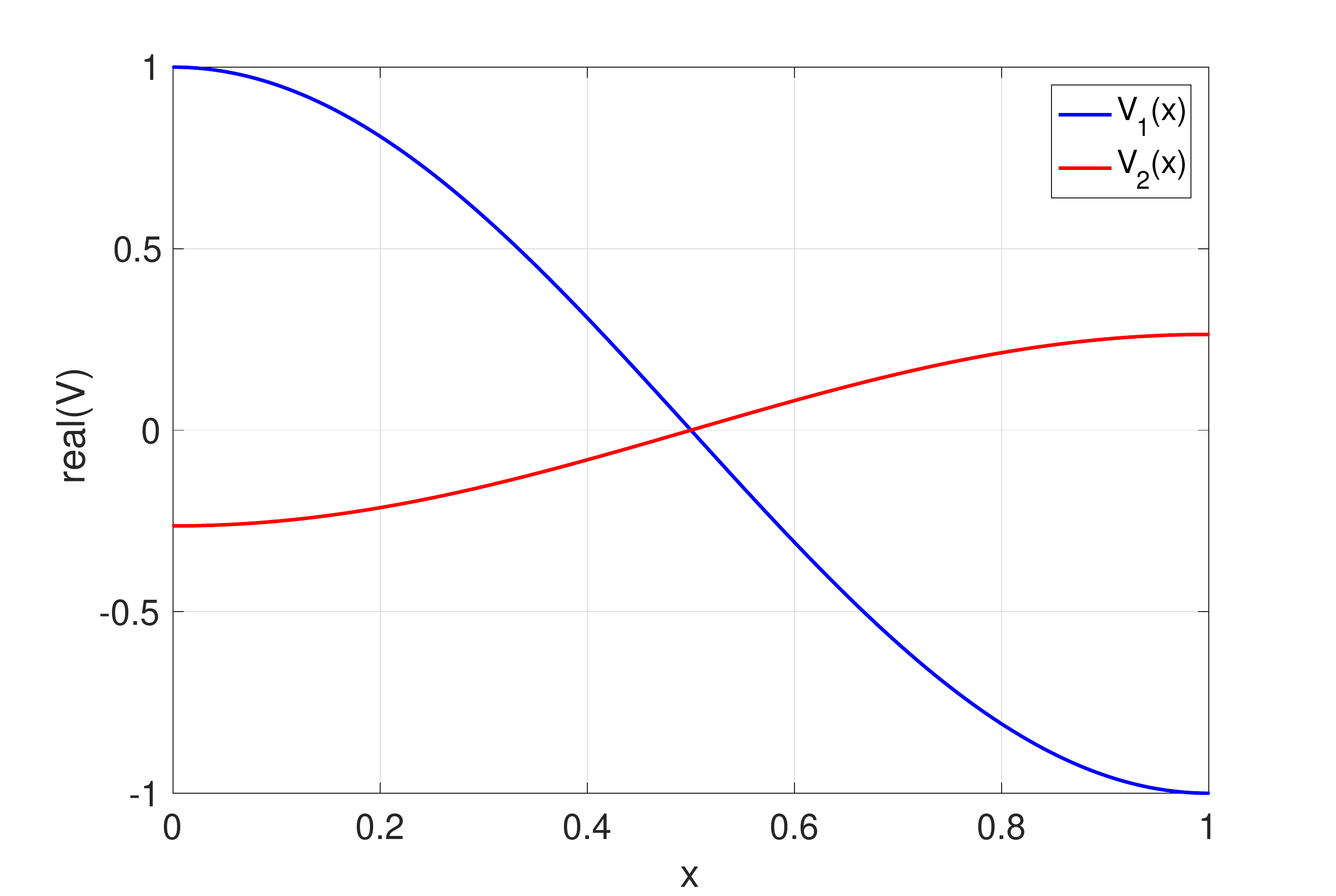}
		\caption{}
		\label{fig:eigfunc}
	\end{subfigure}
	\begin{subfigure}{0.50\textwidth}
		\centering
		\includegraphics[width=\textwidth]{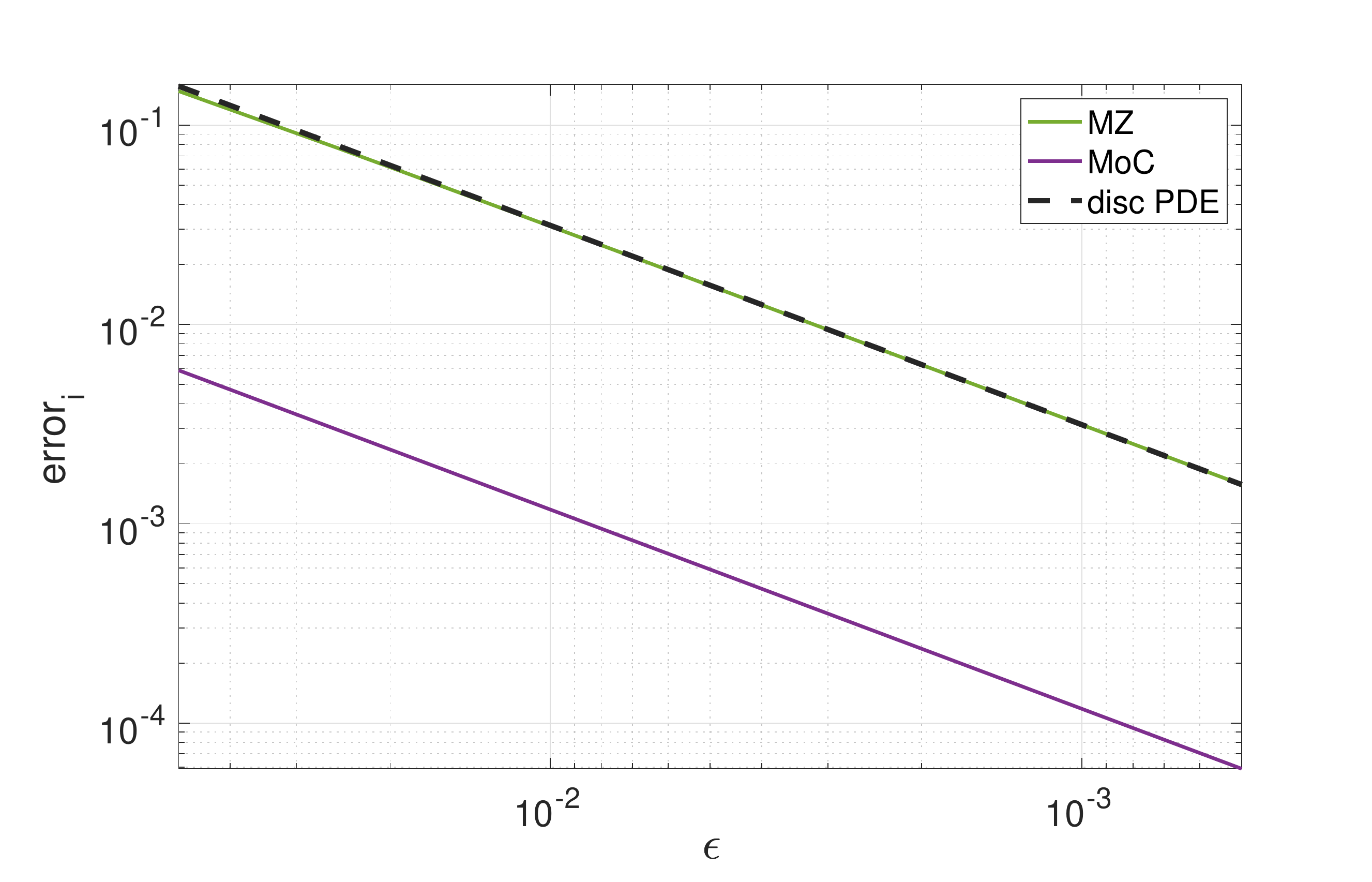}
		\caption{}
		\label{fig:error_BC}
	\end{subfigure}
	\caption{(a) The real part of the eigenfunctions of the PDE and (b) the errors of the discretized PDE (plotted for $\epsilon=\frac{1}{N}$), MoC and MZ model as calculated using the eigenfunctions (with $\alpha=0$).}
	\label{fig:errors}
\end{figure}

Next we compare the performance of the two DDE models derived via the MZ projection (Section \ref{sec:mzamo}\ref{ssec:disc}\ref{ssec:mzproj}\ref{ssec:amonm}) and the MoC (Section \ref{sec:mzamo}\ref{ssec:delchar}) with the exact and discretized PDE models.  We do this through the calculation of the eigenfunctions of the respective models. Figure \ref{fig:eigfunc} shows the real part of the eigenfunctions for each component of the exact PDE calculated using the Chebfun open-source software \cite{Driscoll2014}.  We also compute the eigenfunctions for each approximation to the exact PDE: discretized PDE, delay via method of characteristics (MoC), and delay via Mori-Zwanzig (MZ). We scale all the eigenfunctions such that $V_1(x=0)=1$.  In order to compare the relative approximations, we calculate the error with respect to the eigenfunctions satisfying the PDE boundary conditions \eqref{eq:bcsh}:
\begin{equation} \label{eq:error_BC}
\mathrm{error}_i := |V_1^i(x=0)+V_1^i(x=1)| \qquad \mathrm{for\:} i\in\{\mathrm{disc\:PDE,\:MoC,\:MZ}\}
\end{equation}
As the approximate systems approach the exact PDE, the eigenfunctions of the respective models are expected to converge and therefore satisfy the PDE boundary conditions.  We then would expect \eqref{eq:error_BC} to approach zero as $N$ ($\epsilon$) is increased (decreased) if the models are good approximations of the exact PDE. We calculate \eqref{eq:error_BC} for a range of $N$ and $\epsilon$ values and plot the respective errors in Figure \ref{fig:error_BC}. It can be seen that decreasing the $\epsilon$ term in the delay equations has the same effect on the approximation to the exact PDE as increasing $N$ in the discretization. In other words, the error introduced through the $\epsilon$ term in the `smoothing' approximation of delay difference equations is proportional to the error introduced by discretization methods of wave equation PDEs ($\epsilon \propto 1/N$).

\section{Summary, Discussion and Conclusion}
\label{sec:disc}

A delay model for the AMO has been derived from a three-layer model by S\'evellec and Huck \cite{Sevellec2015} using the MZ formalism. This formalism gives a rewriting of a system of ordinary differential equations \cite{Chorin2002} which contains a Markovian, a noise and a memory term. The advantage of this delay model, over for example mode-decomposition, is that it precisely shows the propagating wave nature of the AMO, through the (inverse) travel times $l_\pm$, hence providing more support for the  thermal Rossby wave mechanism proposed in \cite{TeRaa2002}. In a similar way as the ENSO delayed oscillator model \cite{Suarez1988, Battisti1989}, this model can also be used to study effects of background state, non-stationary forcing, noise, and possibly state-dependent delay versions, on the behavior of the AMO \cite{Tziperman1994,Keane2019}.

The derived model for the AMO is a first order delay difference model, in contrast to the delay differential model for ENSO \cite{Falkena2019}. This means that the current state is fully determined by past states. This type of model can exhibit an increasing switching frequency between states \cite{Ghil2008}, making it physically unrealistic. Hence, an $\epsilon \frac{\d}{\d t}$-term was added to prevent this behaviour and allow for better numerical treatment. We were able to derive an error term for this approximation using the MZ formalism, and relate this error to the upwind discretization scheme used in solving the original PDE model. For the non-damped version of the AMO model ($\alpha=0$) the MZ formalism is not strictly necessary for deriving a delay model, although it is more general and can be extended to other types of models. The method of characteristics yields the same delay difference equation, as shown in Section \ref{sec:amodel}\ref{ssec:delchar}. We showed that the way in which the smoothing approximation $\epsilon \frac{\d}{\d t}$-term is added affects the size of the error of the delay-model. Furthermore, the error introduced through the smoothing $\epsilon \frac{\d}{\d t}$-term is proportional to that introduced by discretization methods, as discussed in Section \ref{sec:amodel}\ref{ssec:compdel}.

The PDE model of the AMO by S\'evellec and Huck \cite{Sevellec2015}, the starting point for deriving the delay model, does not contain a background overturning circulation. This results in a high frequency model oscillation which is undesired to study the AMO. As discussed in Section \ref{sec:amo}\ref{ssec:backover}, adding the meridional overturning circulation to the background state of the model results in a damping of this high frequency oscillation.
This also becomes clear when studying the delay model corresponding to this extended AMO model which has been derived in the Supplementary Material. It shows a weakening of the delayed effect of the high-frequency mode (short delay-time), while the low-frequency mode (long delay-time) is enhanced by the extended background state.

The method of deriving delay equations as applied to a PDE model of the  AMO (and also to a PDE  model of ENSO  \cite{Falkena2019}) can be generalized. It is expected that any diagonalizable linear system of wave equations can be rewritten in the form of a delay difference system, given sufficient coupling that allows for the transfer of a signal between either different variables or boundaries. Integration along characteristics would yield the dominant terms, but the MZ formalism additionally gives error terms to a smoothening approximation for solving the delay difference system. The necessity of the diagonalizability remains to be investigated. For non-diagonalizable systems it is not yet clear whether a similar result would hold as the computations become more involved. 

For nonlinear models the method applied here to the AMO model has to be generalized. Although the MZ formalism is general and does not rely on the linearity of the system, the difficulty lies in solving the orthogonal dynamics equation as noted extensively in the literature \cite{Chorin2002, Darve2009, Zhu2018b}. For linear models the pseudo-orthogonal dynamics approximation can be used to arrive at the final result \cite{Gouasmi2017}. For nonlinear models this approximation cannot be used without first showing its accuracy. In \cite{Falkena2019} it was shown to result in a significant error for the nonlinear ENSO model studied. A better approximation for nonlinear models first needs to be proposed before the method described here can be accurately generalized. This step is necessary for a reliable and accurate application of the MZ formalism to nonlinear models of climate phenomena. 

Many PDE models used to describe climate phenomena contain some type of wave dynamics. We have shown in this study that projecting a system of wave equations onto one location yields a delay model. This would imply that more  climate variability phenomena could be described by a delay equation when there is a physical mechanism that suggests memory effects. Once better methods of approximation of the orthogonal dynamics are available, it may be possible to derive accurate nonlinear delay models of climate phenomena, thus clarifying dynamical mechanisms and allowing for further analysis.

\section*{\large{Author Contributions}}
S.K.J.F., C.Q. and H.A.D. designed the study. The work was carried out mainly by S.K.J.F. and C.Q.. J.S. designed the analysis and comparison of the delay models in Section 5. All authors contributed to the work, discussed the results, read and approved the manuscript and agree to be held accountable for the work performed therein.

\section*{\large{Funding}}
This work was supported by funding from the European Union Horizon 2020 research and innovation programme for the ITN CRITICS under Grant Agreement no. 643073 (C.Q., J.S. and H.A.D.). S.K.J.F. was supported by the Centre for Doctoral Training in Mathematics of Planet Earth, UK EPSRC funded (grant EP/L016613/1) and J.S. by EPSRC grants no. EP/N023544/1 and no. EP/N014391/1.




\bibliographystyle{plain}


\clearpage

{\centering \Huge{Supplementary Information}\\[15mm] \par}

\appendix
\beginsupplement

\section{AMO Model}

\subsection{Relation to Overturning Circulation}

The AMO temperature model (3.4) also provides information about oscillations in the zonal and meridional overturning circulation. Using thermal wind balance 
\begin{equation}
\label{eq:thermalwind_ot}
f \partial_z v = \alpha_T g \partial_x T,
\end{equation}
the vertical shear in the meridional flow can be related to the zonal temperature gradient. This indicated there is a quarter phase difference between the two (as seen in \cite{TeRaa2002}). The vertical shear in the meridional flow indicates whether there is a positive or negative perturbation in the meridional overturning. A positive perturbation in $\partial_z v$ means more northward flow at the surface compared to the bottom, and thus a positive perturbation in the overturning circulation. Similarly a negative perturbation in $\partial_z v$ corresponds to a negative overturning perturbation. Note that in the model considered these flow perturbations do not result in temperature perturbations by the assumption of dominant zonal dynamics.

To get an idea of the behaviour of perturbations in the zonal overturning the $y$-averaged continuity equation is considered:
\begin{equation}
\label{eq:continuity_yav}
\partial_x u + \partial_z w = 0,
\end{equation}
where it is assumed that there is no flow through the boundaries of the basin. Taking the $z$-derivative of this equation and using Sverdrup balance
\begin{equation}
\label{eq:sverdrup}
\beta v = f \partial_z w,
\end{equation}
yields the following equation
\begin{equation}
\label{eq:phasediff}
\partial_x \big( \partial_z u \big) = -\frac{\beta}{f} \partial_z v = -\frac{\beta}{f} \frac{\alpha_T g}{f} \partial_x T.
\end{equation}
This indicates there is a difference of half a phase between zonal overturning perturbations and temperature perturbations. It also implies that the phase difference between the zonal and meridional overturning perturbations is a quarter phase, as is expected from literature (\cite{TeRaa2002}).

\begin{figure}[h]
	\centering
	\includegraphics[width=.8\textwidth]{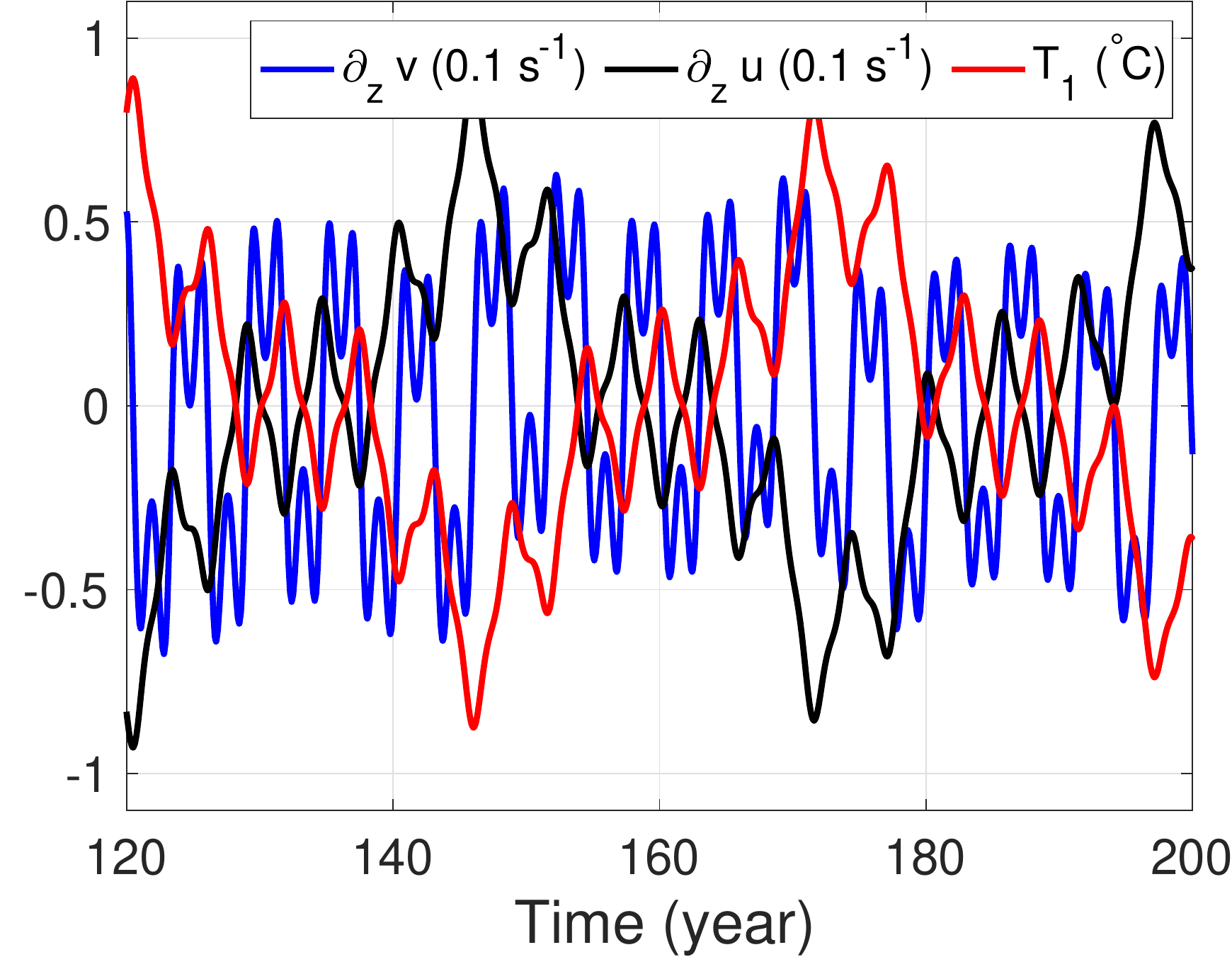}
	\caption{\label{fig:amostr}The evolution of the vertical shear in the meridional and zonal direction (blue) together with the temperature oscillations (red) in the top layer.}
\end{figure}

In Figure \ref{fig:amostr} the evolution of the perturbations in the vertical shear in the zonal and meridional directions is shown as computed from the temperature simulation. For both the short and long period oscillations, the quarter and half phase difference between the temperature oscillations and the meridional and zonal overturning perturbations respectively can be seen. A positive peak in temperature coincides with a negative zonal overturning perturbation, which transports this perturbation westward. Physically a short delay between these two would be expected, but due to the assumption of immediate thermal wind balance this is not present in the model. The resulting zonal temperature gradient causes a negative meridional overturning perturbation with equatorward surface flow lagging by a quarter phase. This is followed by a negative temperature perturbation and a positive zonal overturning, inducing a positive meridional overturning perturbation after which the oscillation starts again. This represents the physical process of the AMO as described in Section 1.

\subsection{Derivation of Extended Background AMO Model}

The temperature equation, including background overturning is: 
\begin{equation}
\partial_t T = -\bar{u} \partial_x T -\bar{v} \partial_y T -\bar{w} \partial_z T - v \partial_y \bar{T} - w \partial_z \bar{T} + \kappa \partial_{xx} T,
\end{equation}
where $\kappa$ is the horizontal eddy diffusivity coefficient.

We start by looking at the $y$-derivative of temperature. $\partial_y T$. The thermal wind balance equations for the meridional and zonal vertical shear are
\begin{equation}
\label{eq:thermalwind}
f \partial_z v = \alpha_T g \partial_x T, \qquad f \partial_z u = - \alpha_T g \partial_y T.
\end{equation}
The second equation means $\partial_y T$ can be replaced by $\partial_z u$ in the temperature equation. This gives
\begin{equation}
\partial_t T = -\bar{u} \partial_x T + \bar{v} \frac{f}{\alpha_T g} \partial_z u -\bar{w} \partial_z T - v \partial_y \bar{T} - w \partial_z \bar{T} + \kappa \partial_{xx} T.
\end{equation}

Now we move or attention to the $z$-derivative of zonal wind $\partial_z u$. We use the Sverdrup balance, given in Equation \eqref{eq:sverdrup} and plug this into the continuity equation. This yields
\begin{equation}
\partial_x u + \partial_y v + \frac{\beta}{f} v = 0.
\end{equation}
Integrating this equation over $y$ and assuming that there is no flow through the northern and southern boundary, so $v|_S=v|_N=0$, gives
\begin{equation}
\partial_x u = - \frac{\beta}{f} v,
\end{equation}
or
\begin{equation}
\partial_x \partial_z u = - \frac{\beta}{f} \partial_z v = - \frac{\beta}{f} \frac{\alpha_T g}{f} \partial_x T.
\end{equation}
This implies that $\partial_z u + \frac{\beta}{f} \frac{\alpha_T g}{f} T$ constant is in $x$. Assuming it is zero initially, the result for the temperature equation is
\begin{equation}
\partial_t T = -\bar{u} \partial_x T - \bar{v} \frac{\beta}{f} T - \bar{w} \partial_z T - v \partial_y \bar{T} - w \partial_z \bar{T} + \kappa \partial_{xx} T.
\end{equation}

The next step is to discretize this system over three layers, as in \cite{Sevellec2015}. For this discretization we follow the same steps as \cite{Sevellec2015}, so we have for $v$
\begin{equation}
\begin{split}
v_1 - v_2 &= \frac{1}{2} (h_1 \partial_z v_1 + h_2 \partial_z v_2), \\
v_2 - v_3 &= \frac{1}{2} (h_2 \partial_z v_2 + h_3 \partial_z v_3). \\
\end{split}
\end{equation}
One more condition is needed to be able to express $v$ in terms of $\partial_z u$ (and thus in terms of $T$) and $\partial_z T$ in terms of $T$. For $v$ we use the baroclinic condition, which in the discretized case reads
\begin{equation}
h_1 v_1 + h_2 v_2 + h_3 v_3 = 0.
\end{equation}
Solving the system of the above three equations for all $v_i$ yields
\begin{equation}
\begin{split}
v_1 &= \frac{1}{2H} ( h_1(h_2+h_3) \partial_z v_1 + h_2(h_2+2h_3 )\partial_z v_2 + h_3^2 \phi_3 ), \\
v_2 &= \frac{1}{2H} ( -h_1^2 \partial_z v_1 + h_2(h_3-h_1) \partial_z v_2 + h_3^2 \partial_z v_3 ), \\
v_3 &= \frac{1}{2H} ( -h_1^2 \partial_z v_1 - h_2(2h_1+h_2) \partial_z v_2 - h_2h_3 \partial_z v_3 ).
\end{split}
\end{equation}
Here thermal wind balance can be used to express $v_i$ in terms of $\partial_x T_i$.

To express the vertical velocity $w$ in terms of temperature we use the Sverdrup balance \eqref{eq:sverdrup} to express it in terms of $v$. With the assumption of zero flow at the bottom, so $w|_{-H}=0$, this yields for the first two layers
\begin{equation}
\begin{split}
w_1 &= -\frac{\beta}{2 f} h_1 v_1, \\
w_2 &= -\frac{\beta}{2 f} ( 2 h_1 v_1 + h_2 v_2). \\
\end{split}
\end{equation}
We assume the third layer is at rest, meaning the value in the third layer is not needed.

For the discretization of temperature over the three layers the equations are:
\begin{equation}
\begin{split}
T_1 - T_2 &= \frac{1}{2} (h_1 \partial_z T_1 + h_2 \partial_z T_2), \\
T_2 - T_3 &= \frac{1}{2} (h_2 \partial_z T_2 + h_3 \partial_z T_3), \\
\end{split}
\end{equation}
where it is assumed that there are no temperature perturbations at the bottom of the basin. We find the vertical gradients of temperature to be
\begin{equation}
\begin{split}
\partial_z T_1 &= \frac{2}{h_1} (T_1 - 2 T_2 + 2 T_3), \\
\partial_z T_2 &= \frac{2}{h_2} (T_2 - 2T_3), \\
\partial_z T_3 &= \frac{2}{h_3} T_3.
\end{split}
\end{equation}

Combining the results above and assuming the third layer is at rest the system we find is
\begin{equation}
\begin{split}
\partial_t T_1 & = a_1 \partial_x T_1 - (\bar{v} \frac{\beta}{f} + \bar{w} \frac{2}{h_1}) T_1 + b_1 \partial_x T_2 + \bar{w} \frac{4}{h_1} T_2 + c_1 \partial_x T_3 - \bar{w} \frac{4}{h_1} T_3 + \kappa \partial_{xx} T_1, \\
\partial_t T_2 & = a_2 \partial_x T_1 + b_2 \partial_x T_2  - (\bar{v} \frac{\beta}{f} + \bar{w}\frac{2}{h_2}) T_2 + c_2 \partial_x T_3 + \bar{w}\frac{4}{h_2} T_3 + \kappa \partial_{xx} T_2, \\
\partial_t T_3 & = \kappa \partial_{xx} T_3,
\end{split}
\end{equation}
with boundary conditions
\begin{equation}
\label{eq:bcsh_si}
T_i|_{West} = - T_i|_{East}, \qquad i=1,2,3.
\end{equation}
The constants here are all positive for physically realistic values and defined by
\begin{equation}
\begin{split}
a_1 &= \frac{\alpha_T g}{2Hf}\Big( -h_1(h_2+h_3) \partial_y \bar{T} + \frac{\beta}{2f} h_1^2(h_2+h_3)\partial_z \bar{T} \Big) - \bar{u}, \\
b_1 &= \frac{\alpha_T g}{2Hf}\Big( -h_2(h_2+2h_3) \partial_y \bar{T} + \frac{\beta}{2f} h_1 h_2(h_2+2h_3)\partial_z \bar{T} \Big), \\
c_1 &= \frac{\alpha_T g}{2Hf}\Big( -h_3^2 \partial_y \bar{T} + \frac{\beta}{2f} h_1 h_3^2 \partial_z \bar{T} \Big), \\
a_2 &= \frac{\alpha_T g}{2Hf}\Big( h_1^2\partial_y \bar{T} + \frac{\beta}{2f} h_1^2(h_2+2h_3)\partial_z \bar{T} \Big), \\
b_2 &= \frac{\alpha_T g}{2Hf}\Big( -h_2(h_3-h_1)\partial_y \bar{T} + \frac{\beta}{2f} (4h_1h_2h_3 + h_2^2 (h_1+h_3))\partial_z \bar{T} \Big) - \bar{u}, \\
c_2 &= \frac{\alpha_T g}{2Hf}\Big( -h_3^2\partial_y \bar{T} + \frac{\beta}{2f} h_3^2(2h_1+h_2)\partial_z \bar{T} \Big). \\
\end{split}
\end{equation}
If, as done for the original system, we neglect the third layer the system becomes
\begin{equation}
\begin{split}
\partial_t T_1 & = a_1 \partial_x T_1 - (\bar{v} \frac{\beta}{f} + \bar{w} \frac{2}{h_1}) T_1 + b_1 \partial_x T_2 + \bar{w} \frac{4}{h_1} T_2 + \kappa \partial_{xx} T_1, \\
\partial_t T_2 & = a_2 \partial_x T_1 + b_2 \partial_x T_2  - (\bar{v} \frac{\beta}{f} + \bar{w}\frac{2}{h_2}) T_2 + \kappa \partial_{xx} T_2.
\end{split}
\end{equation}

\section{Application of the Mori-Zwanzig Formalism}

\subsection{Solution to the Orthogonal Dynamics System}

The solution to the orthogonal dynamics system (4.5) for the AMO model can be written using the eigenvalues and generalized eigenvectors of $M_Q$ (4.7) as given in Equations (4.8) and (4.10). It is
\begin{equation}
\label{eq:amoorthdyn}
\begin{split}
\vec{T}_Q(t) &= e^{\lambda_+ t}\Big(c_+^1\vec{v}_+^1 + c_+^2(t \vec{v}_+^1 + \vec{v}_+^2) + ... + c_+^i \Big(\frac{t^{i-1}}{(i-1)!}\vec{v}_+^1 + \frac{t^{i-2}}{(i-2)!}\vec{v}_+^2 + ...  + \vec{v}_+^i\Big) \\
&\qquad + ... + c_+^{N-1} \Big(\frac{t^{N-2}}{(N-2)!}\vec{v}_+^1 + ... + \vec{v}_+^{N-1}\Big)\Big) \\
&\quad + e^{\lambda_- t}\Big(c_-^1\vec{v}_n^1 + c_-^2(t \vec{v}_-^1 + \vec{v}_-^2) + ... + c_-^i \Big(\frac{t^{i-1}}{(i-1)!}\vec{v}_-^1 + \frac{t^{i-2}}{(i-2)!}\vec{v}_-^2 + ... + \vec{v}_-^i\Big) \\
&\qquad + ... + c_-^{N-1} \Big(\frac{t^{N-2}}{(N-2)!}\vec{v}_-^1 + ... + \vec{v}_-^{N-1}\Big)\Big).
\end{split}
\end{equation}
Here the constants $c_\pm^i$ are determined by the initial conditions. Each generalized eigenvector (4.10) has only components in the directions of $T_1^i$ and $T_2^i$. This means that, to find expressions for the constants, the following system has to be solved for each $i$:
\begin{equation}
\begin{split}
c_+^{i} w_+\Big(\frac{dx}{l_+}\Big)^{i-1} + c_-^{i} w_-\Big(\frac{dx}{l_-}\Big)^{i-1} &= T_1^{i}(0), \\
c_+^{i} \Big(\frac{dx}{l_+}\Big)^{i-1} + c_-^{i} \Big(\frac{dx}{l_-}\Big)^{i-1} &= T_2^{i}(0).
\end{split}
\end{equation}
The solution is
\begin{equation}
\label{eq:amomzconstants}
\begin{split}
c_+^{i} &= \Big(\frac{l_+}{dx}\Big)^{i-1} \cdot \frac{T_1^{i}(0)-w_-T_2^{i}(0)}{w_+-w_-}, \\
c_-^{i} &= -\Big(\frac{l_-}{dx}\Big)^{i-1} \cdot \frac{T_1^{i}(0)-w_+T_2^{i}(0)}{w_+-w_-}.
\end{split}
\end{equation}
This way an analytical solution to the orthogonal dynamics equation for general initial conditions has been found.

\subsection{Limiting Behaviour}

\begin{figure}[h]
	\centering
	\includegraphics[width=.7\textwidth]{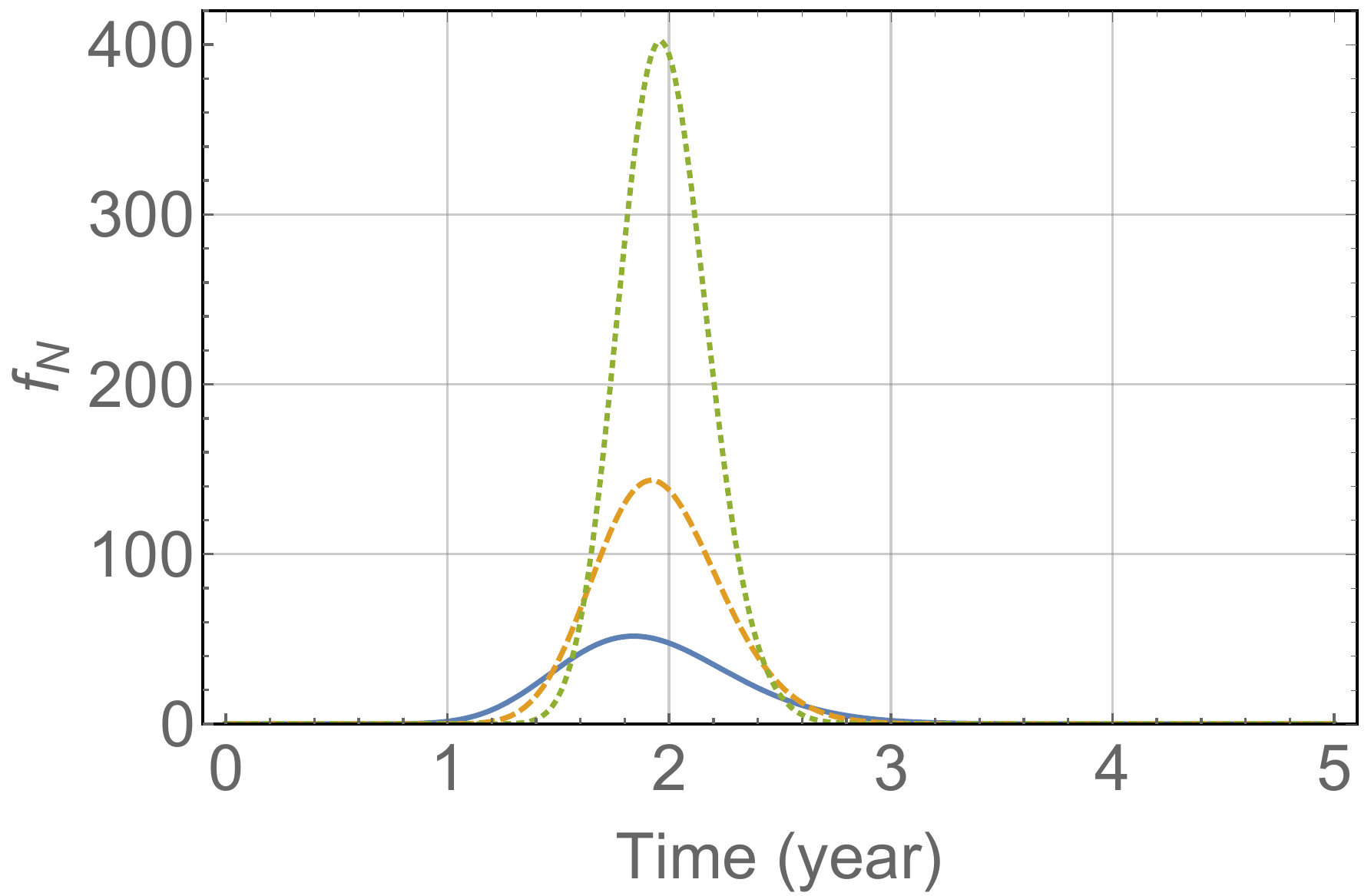}
	\caption{\label{fig:amokern}The function $f_N(t)$ with $\mu=0.5$ for $N=25$ (blue, solid), $N=50$ (yellow, dashed) and $N=100$ (green, dotted).}
\end{figure}

In studying the limiting behaviour of Equation (4.18) we focus on the memory integral first, as this is the term that is expected to result in a delay (\cite{Falkena2019}). The dependence on $N$ of all terms in the memory kernel can be described by the function
\begin{equation}
f_N(t) = N^2 \frac{t^{N-2}}{(N-2)!} e^{-\alpha t} (\mu N)^{N-2} e^{-\mu N t},
\end{equation}
for $\mu=l_\pm$. In Figure \ref{fig:amokern} this function is plotted for several $N$ and fixed $\mu$ showing a peak at $1/\mu$ which increases with height as $N$ gets larger. This behaviour is caused by the projection onto one location, which prevents waves from travelling through the basin in the orthogonal dynamics system. The result is an accumulation of energy at the location of the resolved variables. Because one integrates over the memory kernel to obtain the memory term the blow up in peak height does not necessarily imply a blow up of the memory term itself.

To find an approximate expression for the memory integral we use Laplace's approximation (\cite{Erdelyi1956}), which states that for a large number $N$, a smooth function $h(x)$, and a twice differentiable function $g(x)$, the following integral can be approximated by
\begin{equation}
\label{eq:laplace}
\begin{split}
\int_a^b h(x) e^{Ng(x)} \d x &= \sqrt{\frac{2\pi}{N|g''(x_0)|}} e^{Ng(x_0)} \cdot \Big( h(x_0) + \frac{1}{N}\Big( -\frac{h''(x_0)}{2g''(x_0)} + \frac{h(x_0)g''''(x_0)}{8(g''(x_0))^2} \\
&\qquad + \frac{h'(x_0)g'''(x_0)}{2(g''(x_0))^2} - \frac{5h(x_0)(g'''(x_0))^2}{24(g''(x_0))^3}\Big) + \mathcal{O}\Big(\frac{1}{N^2}\Big)\Big),
\end{split}
\end{equation}
provided that there is an $x_0\in(a,b)$ such that $g(x)$ is only close to $g(x_0)$ if $x$ is close to $x_0$. At that point $g(x)$ is required to have a maximum, i.e. $g''(x_0)<0$.

All components in the memory integrals (4.18) can be written in the form $h_{ij\pm}(s)e^{Ng_{\pm}(s)}$, for $i=1,2$, $j=1,2$, with
\begin{equation}
\begin{split}
g_\pm(s) &= \log (l_\pm (t-s)) - l_\pm(t-s), \\
h_{ij\pm}(s) &= \frac{N^N}{(N-2)!} (l_\pm (t-s))^{-2} e^{-\alpha(t-s)} C_{j\pm} T_i(s),
\end{split}
\end{equation}
where $C_{j\pm}$ is either $A_{j\pm}$ or $B_{j\pm}$ as given in Equation (4.17). The function $h_{ij\pm}$ is smooth meaning the first of the conditions for Laplace's approximation is satisfied. To verify the other two conditions the first two derivatives of $g_\pm(s)$ are required. We find the first derivative exists and is only zero for $s_0=t-\frac{1}{l_\pm}$ indicating there is only one extreme. This means that $g_\pm(s)$ is only close to $g_\pm(s_0)$ for $s$ near $s_0$, satisfying the second condition. The second derivative is $g_\pm''(s_0) = -l_\pm^2 <0$, indicating this extreme is a maximum meaning the third condition is satisfied. Therefore we are justified in applying Laplace's approximation, which yields:
\begin{equation}
\label{eq:Ndepmemory}
\begin{split}
\int_0^t h_{ij\pm}(s) e^{Ng_\pm(s)} \d s &= \sqrt{\frac{2\pi}{N l_\pm^2}}e^{-N} \frac{N^N}{(N-2)!} e^{-\frac{\alpha}{l_\pm}} C_{j\pm} \Big( T_i\big( t-\frac{1}{l_\pm} \big) \\
&\qquad + \frac{1}{N} \frac{1}{2l_\pm^2} \Big( (2 l_\pm(l_\pm+\alpha)+\alpha^2) T_i\big( t-\frac{1}{l_\pm} \big) + 2(l_\pm+\alpha) T_i'\big( t-\frac{1}{l_\pm} \big) \\
&\qquad + T_i''\big( t-\frac{1}{l_\pm} \big) \Big) + \mathcal{O}\Big(\frac{1}{N^2}\Big) \Big).
\end{split}
\end{equation}
In this equation delay terms have emerged due to the location of the maximum of $g_\pm(s)$ at $s_0=t-\frac{1}{l_\pm}$. This means that the AMO can be modelled by some type of delay equation, similar to the results found for ENSO by \cite{Falkena2019}.

To simplify this expression we expand the first part of Equation \eqref{eq:Ndepmemory} using a Taylor expansion. Expanding in $1/N$ around zero we find
\begin{equation}
\sqrt{\frac{2\pi}{N l_\pm^2}}e^{-N} \frac{N^N}{(N-2)!} = \frac{N}{l_\pm} \Big( 1 - \frac{13}{12}\frac{1}{N} + \mathcal{O}\Big( \frac{1}{N^2} \Big) \Big),
\end{equation}
allowing for the expansion of the memory integrals in terms of order $\epsilon = \frac{1}{N}$. Equation \eqref{eq:Ndepmemory} thus becomes
\begin{equation}
\begin{split}
\int_0^t h_{ij\pm}(s) e^{Ng_\pm(s)} \d s &= \frac{N}{l_\pm} e^{-\frac{\alpha}{l_\pm}} C_{j\pm} \Big( T_i\big( t-\frac{1}{l_\pm} \big) + \epsilon \frac{1}{2l_\pm^2} \Big( \big((l_\pm+\alpha)^2-\frac{7}{6}l_\pm^2\big) T_i\big( t-\frac{1}{l_\pm} \big) \\
&\qquad + 2(l_\pm+\alpha) T_i'\big( t-\frac{1}{l_\pm} \big) + T_i''\big( t-\frac{1}{l_\pm} \big) \Big) + \mathcal{O}(\epsilon^2) \Big).
\end{split}
\end{equation}
Here we keep one factor of $N=\frac{1}{\epsilon}$ to align with the Markovian term in Equation (4.18). We now have concise expressions for the memory terms in Equation (4.18).

The behaviour of the noise terms for increasing $N$ also needs to be investigated. Each component of the sum in the noise terms is proportional to
\begin{equation}
\label{eq:noiseNdep}
N e^{-l_\pm N t} \frac{(l_\pm N t)^{N-k}}{(N-k)!},
\end{equation}
for some $k=2,...,N$. Dividing by $N$ the resulting function peaks for $t=\frac{1}{l_\pm}$, just as the memory kernels. At this $t$ the expression in Equation \eqref{eq:noiseNdep} can be expanded in terms of order $\epsilon$, similarly to what has been done for the memory terms:
\begin{equation}
N e^{-N} \frac{(N)^{N-k}}{(N-k)!} = \frac{1}{\sqrt{2\pi}}\frac{N}{\sqrt{\epsilon}}( 1-\frac{13}{12}\epsilon + \mathcal{O}(\epsilon^2) ).
\end{equation}
We find that the noise terms decay an order $\sqrt{\epsilon}$ faster than the memory terms (at first order). Note that this is an effect occurring at one time only as a remnant of the initial conditions. At all other times the noise term does not have any effect since Equation \eqref{eq:noiseNdep} approaches zero faster than $\epsilon^2$ for all other times. Therefore we neglect the noise term in the final equations of applying the MZ formalism to the AMO model.

To arrive at the final reduced model for the AMO we divide Equation (4.18) by $N$ and use Equation \eqref{eq:Ndepmemory} to write down expressions for the memory terms. We arrive at
\begin{equation}
\label{eq:Amzres}
\begin{split}
\epsilon \frac{\d T_1}{\d t} &= -a_1 T_1(t) - b_1 T_2(t) + \frac{A_{1+}}{l_+} e^{-\frac{\alpha}{l_+}} T_1\big( t-\frac{1}{l_+} \big) + \frac{B_{1+}}{l_+} e^{-\frac{\alpha}{l_+}} T_2\big( t-\frac{1}{l_+} \big) \\
&\quad + \frac{A_{1-}}{l_-} e^{-\frac{\alpha}{l_-}} T_1\big( t-\frac{1}{l_-} \big) + \frac{B_{1-}}{l_-} e^{-\frac{\alpha}{l_-}} T_2\big( t-\frac{1}{l_-} \big) + \epsilon f_{\epsilon1}(t) + \mathcal{O}(\epsilon^2), \\
\epsilon \frac{\d T_2}{\d t} &= -a_2 T_1(t) - b_2 T_2(t) + \frac{A_{2+}}{l_+} e^{-\frac{\alpha}{l_+}} T_1\big( t-\frac{1}{l_+} \big) + \frac{B_{2+}}{l_+} e^{-\frac{\alpha}{l_+}} T_2\big( t-\frac{1}{l_+} \big) \\
&\quad + \frac{A_{2-}}{l_-} e^{-\frac{\alpha}{l_-}} T_1\big( t-\frac{1}{l_-} \big) + \frac{B_{2-}}{l_-} e^{-\frac{\alpha}{l_-}} T_2\big( t-\frac{1}{l_-} \big) + \epsilon f_{\epsilon2}(t) + \mathcal{O}(\epsilon^2),
\end{split}
\end{equation}
where
\begin{equation}
\label{eq:Amzerror}
\begin{split}
f_{\epsilon1}(t) &=  - \alpha T_1(t) +\frac{A_{1+}}{l_+} e^{-\frac{\alpha}{l_+}}g_{\epsilon+}\big(T_1\big) +\frac{B_{1+}}{l_+} e^{-\frac{\alpha}{l_+}} g_{\epsilon+}\big(T_2\big) \\
&\qquad +\frac{A_{1-}}{l_-} e^{-\frac{\alpha}{l_-}} g_{\epsilon-}\big(T_1\big)+ \frac{B_{1-}}{l_-} e^{-\frac{\alpha}{l_-}} g_{\epsilon-}\big(T_2\big), \\
f_{\epsilon2}(t) &= - \alpha T_2(t) +\frac{A_{2+}}{l_+} e^{-\frac{\alpha}{l_+}} g_{\epsilon+}\big(T_1\big) +\frac{B_{2+}}{l_+} e^{-\frac{\alpha}{l_+}} g_{\epsilon+}\big(T_2\big)\\
&\qquad +\frac{A_{2-}}{l_-} e^{-\frac{\alpha}{l_-}} g_{\epsilon-}\big(T_1\big) + \frac{B_{2-}}{l_-} e^{-\frac{\alpha}{l_-}} g_{\epsilon-}\big(T_2\big),
\end{split}
\end{equation}
for
\begin{equation}
g_{\epsilon\pm}\big(T\big) = \frac{1}{2l_\pm^2} \Big( \big((l_\pm+\alpha)^2-\frac{7}{6}l_\pm^2\big) T\big( t-\frac{1}{l_\pm} \big) + 2(l_\pm+\alpha) T'\big( t-\frac{1}{l_\pm} \big) + T''\big( t-\frac{1}{l_\pm} \big) \Big).
\end{equation}

\section{Application of the Mori-Zwanzig Formalism to Extended AMO Model}

In this section of the supplementary material we go through the derivation of a delay model for the extended AMO model as described in Section 3(b). This derivation goes through the same steps as Section 4, but is more involved in that the terms are lengthier. We start by discretizing the model, projecting onto $T_1^0$ and $T_2^0$ and then turn to the computation of the memory terms by solving the orthogonal dynamics equation. Note that we do not discuss the noise terms, as we do not need it for the final delay model. Its computation is equivalent to the one in the paper, but know using the new eigenvalues and eigenvectors of the orthogonal dynamics system.

\subsection{Discretization}

The discretized equations for the extended AMO model are
\begin{equation}
\label{eq:amo_disc}
\begin{split}
\partial_t T_1^n &= \frac{a_1}{dx}(T_1^{n+1}-T_1^n) + \frac{b_1}{dx}(T_2^{n+1}-T_2^n) - (\alpha+\beta_1) T_1^n - \beta_2 T_2^n, \\
\partial_t T_2^n &= \frac{a_2}{dx}(T_1^{n+1}-T_1^n) + \frac{b_2}{dx}(T_2^{n+1}-T_2^n) - (\alpha+\beta_3) T_2^n,
\end{split}
\end{equation}
for $n=0,...,N$ (such that $T_i^k\approx T_i(k/N)$ and $dx=1/N$), with boundary conditions
\begin{equation}
T_1^{N} = -T_1^0,  \qquad T_2^{N} = -T_2^0.
\end{equation}
The main difference with the model without background overturning is the additional coupling between the two layers via the $\beta_2 T_2^n$ term in the equation for $T_1^n$. The eigenvalues show the same qualitative behaviour as for the AMO model without background overturning circulation relative to the eigenvalues of the PDE. Note here that for $\alpha=0$ the model with overturning is slightly unstable (eigenvalue with positive real-part) for the low-frequency modes, while the high-frequency modes are damped more strongly.

\subsection{Markovian Terms}

Following the results for the standard AMO model, we project again onto $T_1^0$ and $T_2^0$. The eigenvalues of the orthogonal dynamics are discussed in the following section. Here we give the Markovian term, which reads:
\begin{equation}
\begin{split}
R_{T_1^0} &= -\frac{a_1}{dx}T_1^0 - \frac{b_1}{dx}T_2^0 - (\alpha+\beta_1) T_1^0 - \beta_2 T_2^0, \\
R_{T_2^0} &= -\frac{a_2}{dx}T_1^0 - \frac{b_2}{dx}T_2^0 - (\alpha+\beta_3) T_2^0.
\end{split}
\end{equation}
This contains some additional terms compared to the standard AMO model. These terms will be of order $\epsilon$ in the final delay model, as they are similar to the $\alpha$-terms in the original AMO model.

\subsection{Memory Term}

With the projection onto $T_1^0$ and $T_2^0$ the orthogonal dynamics system reads
\begin{equation}
\partial_t \vec{T}_Q = M_Q \vec{T}_Q,
\end{equation}
with of the same form $M_Q$ as to the matrix discussed in the main paper, but here the blocks on the diagonal read
\begin{equation}
\begin{pmatrix}
-\frac{a_1}{dx}-\alpha-\beta_1 & -\frac{b_1}{dx}-\beta_2 \\
-\frac{a_2}{dx} & -\frac{b_2}{dx}-\alpha-\beta_3
\end{pmatrix}.
\end{equation}
The eigenvalues of this matrix are
\begin{equation}
\begin{split}
\lambda_\pm &= -\alpha - \frac{1}{2} \Big( \frac{a_1+b_2}{dx} + \beta_1+\beta_3 \\
&\qquad \qquad \qquad \pm \sqrt{ \frac{(a_1-b_2)^2 + 4a_2b_1}{dx^2} + \frac{2(a_1-b_2)(\beta_1-\beta_3) + 4a_2\beta_2}{dx} + (\beta_1-\beta_3)^2} \quad \Big),
\end{split}
\end{equation}
which as for the standard AMO model can be written in the form
\begin{equation}
\lambda_\pm = -\alpha - \frac{l_\pm}{dx},
\end{equation}
with
\begin{equation}
\begin{split}
l_\pm &= \frac{1}{2} \Big( a_1+b_2 + dx(\beta_1+\beta_3) \\
&\qquad \qquad \pm \sqrt{ (a_1-b_2)^2 + 4a_2b_1 + 2dx ((a_1-b_2)(\beta_1-\beta_3) + 2a_2\beta_2) + dx^2(\beta_1-\beta_3)^2} \Big).
\end{split}
\end{equation}
Note that for $\beta_{1,2,3}=0$ we obtain exactly the same expression as for the standard AMO model, and the dependence of $l_\pm$ on $dx$ disappears. As for the AMO model without background overturning circulation we have the eigenvectors for $i = 1,...,N-1$
\begin{equation}
\vec{v}_\pm^i = \Big( \frac{dx}{l_\pm} \Big)^{i-1} \cdot (0,...,0,w_\pm,1,0,...,0),
\end{equation}
where the non-zero values are located on the coordinates corresponding to location $i$. Here for the extended AMO model
\begin{equation}
\begin{split}
w_\pm &= \frac{1}{2a_2}\Big(a_1-b_2+ dx(\beta_1-\beta_3) \\
&\qquad \qquad \pm \sqrt{(a_1-b_2)^2 + 4a_2b_1 + 2dx ((a_1-b_2)(\beta_1-\beta_3) + 2a_2\beta_2) + dx^2(\beta_1-\beta_3)^2}\Big).
\end{split}
\end{equation}

The expressions for the noise and memory terms are the same as for the standard AMO model (Equations (4.13), (4.16), (4.18)), but with $l_\pm$ and $w_\pm$ as defined above. The dependence of both $l_\pm$ and $w_\pm$ on $dx$, or $N=\frac{1}{dx}$, complicates the computation of the limit behaviour as $N\rightarrow\infty$. To consider this limit we first use a Taylor expansion to expand $\sqrt{A+Bx+Cx^2}$ and $(\sqrt{A+Bx+Cx^2})^{-1}$ around $x=0$ in terms of order $x$. The second expansion is needed due to the $\frac{1}{w_+-w_-}$ terms in the memory integrand (Equation (4.17)). We have
\begin{equation}
\begin{split}
(A+Bx+Cx^2)^{1/2} &\approx A^{1/2} + \frac{1}{2} A^{-1/2}B x + \frac{1}{2}(-\frac{1}{4} A^{-3/2}B^2 + A^{-1/2}C)x^2 \\
&\qquad + \frac{1}{6}(\frac{3}{8} A^{-5/2}B^3 -\frac{3}{2} A^{-3/2}BC)x^3 + ...,
\end{split}
\end{equation}
and
\begin{equation}
\begin{split}
(A+Bx+Cx^2)^{-1/2} &\approx A^{-1/2} - \frac{1}{2} A^{-3/2}B x + \frac{1}{2}(\frac{3}{4} A^{-5/2}B^2 - A^{-3/2}C)x^2 \\
&\qquad + \frac{1}{6}(-\frac{15}{8} A^{-7/2}B^3 +\frac{9}{2} A^{-5/2}BC)x^3 + ...,
\end{split}
\end{equation}
where
\begin{equation}
\begin{split}
A &= (a_1-b_2)^2 + 4a_2b_1, \\
B &= 2((a_1-b_2)(\beta_1-\beta_3) + 2a_2\beta_2), \\
C &= (\beta_1-\beta_3)^2.
\end{split}
\end{equation}
This allows us to write
\begin{equation}
\begin{split}
l_\pm &= l_\pm^0 + dx \cdot l_\pm^1 + dx^2 \cdot l_\pm^2 + dx^3 \cdot l_\pm^3 + ..., \\
w_\pm &= w_\pm^0 + dx \cdot w_\pm^1 + dx^2 \cdot w_\pm^2 + dx^3 \cdot w_\pm^3 + ..., \\
\frac{1}{w_+-w_-} &= u^0 + dx \cdot u^1 + dx^2 \cdot u^2 + dx^3 \cdot u^3 + ...,
\end{split}
\end{equation}
where
\begin{equation}
\begin{split}
l_\pm^0 &= \frac{1}{2} (a_1+b_2 \pm \sqrt{(a_1-b_2)^2 + 4a_2b_1}), \\
l_\pm^1 &= \frac{1}{2} \Big(\beta_1+\beta_3 \pm \frac{(a_1-b_2)(\beta_1-\beta_3)+2a_2\beta_2}{\sqrt{(a_1-b_2)^2 + 4a_2b_1}}\Big), \\
l_\pm^2 &= \pm a_2 \frac{b_1(\beta_1-\beta_3)^2-\beta_2(a_1-b_2)(\beta_1-\beta_3)-a_2\beta_2^2}{((a_1-b_2)^2 + 4a_2b_1)^{3/2}}, \\
l_\pm^3 &= \pm \frac{(a_1-b_2)(\beta_1-\beta_3)+2a_2\beta_2}{32((a_1-b_2)^2 + 4a_2b_1)^{5/2}}\cdot\Big(\big((a_1-b_2)(\beta_1-\beta_3)+2a_2\beta_2\big)^2 \\
&\qquad \qquad \qquad \qquad \qquad \qquad \qquad \qquad -4\big(\beta_1-\beta_3)^2((a_1-b_2)^2 + 4a_2b_1\big)\Big), \\
w_\pm^0 &= \frac{1}{2a_2} (a_1-b_2 \pm \sqrt{(a_1-b_2)^2 + 4a_2b_1}), \\
w_\pm^1 &= \frac{1}{2a_2} \Big(\beta_1-\beta_3 \pm \frac{(a_1-b_2)(\beta_1-\beta_3)+2a_2\beta_2}{\sqrt{(a_1-b_2)^2 + 4a_2b_1}}\Big), \\
w_\pm^{2,3} &= \frac{l_\pm^{2,3}}{a_2}, \\
u^0 &= \frac{a_2}{\sqrt{(a_1-b_2)^2 + 4a_2b_1}}, \\
u^1 &= a_2 \frac{-(a_1-b_2)(\beta_1-\beta_3)-2a_2\beta_2}{((a_1-b_2)^2 + 4a_2b_1)^{3/2}}, \\
u^2 &= \frac{1}{8}\frac{3((a_1-b_2)(\beta_1-\beta_3)+2a_2\beta_2)^2 - 4(\beta_1-\beta_3)^2((a_1-b_2)^2 + 4a_2b_1)}{((a_1-b_2)^2 + 4a_2b_1)^{5/2}}, \\
u^3 &= \frac{1}{48}\Big( \frac{-15((a_1-b_2)(\beta_1-\beta_3)+2a_2\beta_2)^3}{((a_1-b_2)^2 + 4a_2b_1)^{7/2}} + \frac{32(\beta_1-\beta_3)^2((a_1-b_2)(\beta_1-\beta_3)+2a_2\beta_2)}{((a_1-b_2)^2 + 4a_2b_1)^{5/2}} \Big).
\end{split}
\end{equation}
For simplicity we take out the first order term $l^1_\pm$ out of the expansion of $l_\pm$ and add it to $\alpha$ in the definition of $\lambda_\pm$, that is:
\begin{equation}
\begin{split}
\lambda_\pm &= -\alpha - l^1_\pm - \frac{l_\pm}{dx}, \\
l_\pm &= l_\pm^0 + dx^2 \cdot l_\pm^2 + dx^3 \cdot l_\pm^3 + ... \quad .
\end{split}
\end{equation}
From here on we use these updated definitions of $\lambda_\pm$ and $l_\pm$.

All terms in the memory integrand are of a similar form, but with different constants. This general form with $\mu=l_\pm$ is
\begin{equation}
\begin{split}
K &\sim N^2 \frac{t^{N-2}}{(N-2)!} e^{-(\alpha+l^1_\pm) t} (\mu^0 N + \mu^2 \frac{1}{N} + \mu^3 \frac{1}{N^2} + ...)^{N-2} \cdot e^{-(\mu^0 N + \mu^2 \frac{1}{N} + \mu^3 \frac{1}{N^2} + ...)t} \\
& \qquad \cdot (c_1(w_+^0 + w_+^1\frac{1}{N} + w_+^2\frac{1}{N^2} + ...)+c_2)(c_3(w_-^0 + w_-^1\frac{1}{N} + w_-^2\frac{1}{N^2} + ...)+c_4) \\
&\qquad \cdot (u^0 + u^1\frac{1}{N} + u^2\frac{1}{N^2} + ...), \\
&\sim N^2 \frac{t^{N-2}}{(N-2)!} e^{-(\alpha+l^1_\pm) t} (\mu^0 N + \mu^2 \frac{1}{N} + \mu^3 \frac{1}{N^2} + ...)^{N-2} \cdot e^{-\mu^0 Nt} \cdot e^{-\mu^2 \frac{1}{N}t} \cdot ... \\
&\qquad \cdot (v^0 + v^1\frac{1}{N} + v^2\frac{1}{N^2} + ...),
\end{split}
\end{equation}
where the $v^i$ are obtained by expanding the product of the final three terms. We note that $v^0$ yields the result for the standard AMO model. As $N\rightarrow\infty$ all exponential terms with $\mu^i$ where $i>1$ will go to one and can be neglected at zeroth and first order, which can be be verified by retaining them when using Laplace's approximation to the memory integral and looking at the leading order terms. This leaves us with
\begin{equation}
K \sim N^2 \frac{t^{N-2}}{(N-2)!} e^{-(\alpha+\mu^1) t} (\mu^0 N + \mu^2 \frac{1}{N} + ...)^{N-2} \cdot e^{-\mu^0 Nt} \cdot (v^0 + v^1\frac{1}{N} + v^2\frac{1}{N^2} + ...).
\end{equation}
The main difficulty in this expression compared to taking the limit for the standard AMO model is the fourth term, being a sum to the power $N-2$. We can write this as
\begin{equation}
\begin{split}
(\mu^0 N + \mu^2 \frac{1}{N} + ...)^{N-2} &= (\mu^0 N)^{N-2} + {{N-2}\choose{1}}(\mu^0 N)^{N-3}(\mu^2 \frac{1}{N} + ...) \\
&\qquad + {{N-2}\choose{2}}(\mu^0 N)^{N-4}(\mu^2 \frac{1}{N} + ...)^2 + ... \quad .
\end{split}
\end{equation}
The leading term $(\mu^0 N)^{N-2}$ is of order $\mathcal{O}(N^{N-2})$. Similarly one can look at the $\mathcal{O}(N^{N-3})$ terms. Separating out the $\mathcal{O}(N^{N-2})$ and $\mathcal{O}(N^{N-3})$ terms yields
\begin{equation}
\begin{split}
(\mu^0 N &+ \mu^2 \frac{1}{N} + ...)^{N-2} &= (\mu^0 N)^{N-2} + (\mu^0 N)^{N-3}\mu^2 + \mathcal{O}\big((\mu^0 N)^{N-4}\big). \\
\end{split}
\end{equation}
This yields the memory integrand:
\begin{equation}
\begin{split}
K &\sim N^2 \frac{t^{N-2}}{(N-2)!} e^{-(\alpha+\mu^1) t} e^{-\mu^0 Nt} (\mu^0 N)^{N-2} v^0 + \frac{1}{N} N^2 \frac{t^{N-2}}{(N-2)!} e^{-(\alpha+\mu^1) t} e^{-\mu^0 Nt} (\mu^0 N)^{N-2} v^1 \\
&\qquad + N^2 \frac{t^{N-2}}{(N-2)!} e^{-(\alpha+\mu^1) t} e^{-\mu^0 Nt} (\mu^0 N)^{N-3} \mu^2 v^0 + \mathcal{O}(\frac{1}{N^2}).
\end{split}
\end{equation}
The leading term is very similar to the one for the standard AMO model, only now $e^{-\alpha t}$ has been replaced by $e^{-(\alpha+\mu^1) t}$. Using Laplace's approximation on the leading order term as before yields up to first order
\begin{equation}
\begin{split}
\int_0^t K(t-s) T_i(s) ds &= \sqrt{\frac{2\pi}{N(\mu^0)^2}} e^{-N} \frac{N^N}{(N-2)!} e^{-\frac{\alpha+\mu^1}{\mu^0}} v^0 \Big( T_i(t-\frac{1}{\mu^0}) \\
&\qquad + \frac{1}{N} \frac{1}{2(\mu^0)^2} \big( (2\mu^0(\mu^0+\alpha+\mu^1) + (\alpha+\mu^1)^2)T_i(t-\frac{1}{\mu^0}) \\
&\qquad \qquad \qquad \qquad  + 2(\mu^0+\alpha+\mu^1)T_i'(t-\frac{1}{\mu^0}) + T_i''(t-\frac{1}{\mu^0}) \big) \Big) \\
&\qquad + \mathcal{O}(\frac{1}{N^2}).
\end{split}
\end{equation}
Similarly, we can compute the additional first order terms by applying Laplace's approximation to the first order terms in the memory integrand:
\begin{equation}
\begin{split}
\frac{1}{N} \int_0^t K(t-s) T_i(s) ds &= \frac{1}{N} \sqrt{\frac{2\pi}{N(\mu^0)^2}} e^{-N} \frac{N^N}{(N-2)!} e^{-\frac{\alpha+\mu^1}{\mu^0}} v^1 T_i(t-\frac{1}{\mu^0}) \\
&\qquad + \frac{1}{N} \sqrt{\frac{2\pi}{N(\mu^0)^3}} e^{-N} \frac{N^N}{(N-2)!} e^{-\frac{\alpha+\mu^1}{\mu^0}} \mu^2 v^0 T_i(t-\frac{1}{\mu^0}) + \mathcal{O}(\frac{1}{N^2}).
\end{split}
\end{equation}
Thus we now can write down the final equation for the memory term in terms of order $\epsilon=\frac{1}{N}$ as has been done before for the standard AMO model. We find for the memory term:
\begin{equation}
\begin{split}
\int_0^t K(t-s) T_i(s) ds &= \frac{N}{\mu^0} e^{-\frac{\alpha+\mu^1}{\mu^0}} \Big( v^0 T_i(t-\frac{1}{\mu^0}) \\
&\qquad +\epsilon \frac{v^0}{2(\mu^0)^2} \big( ((\mu^0+\alpha+\mu^1)^2 - \frac{7}{6}(\mu^0)^2)T_i(t-\frac{1}{\mu^0}) \\
&\qquad \qquad \qquad \qquad  + 2(\mu^0+\alpha+\mu^1)T_i'(t-\frac{1}{\mu^0}) + T_i''(t-\frac{1}{\mu^0}) \big) \\
&\qquad + \epsilon \frac{1}{\mu^0} \big( v^1  + \frac{v^0 \mu^2}{\sqrt{\mu^0}} \big) T_i(t-\frac{1}{\mu^0}) \Big) + \mathcal{O}(\frac{1}{N^2}).
\end{split}
\end{equation}
This also allows us to write down the full delay equations for the AMO model including background overturning circulation. This is equivalent to the approach for the original AMO model and is left to the reader.

We now turn to discussing the effect of the addition of $l^1_\pm$ in the leading order term. We have that $l^1_+$, which corresponds to the short delay, is positive thus leading to an additional damping of this high-frequency mode. On the other hand $l^1_-$, corresponding to the long delay, is negative. This means that without damping $\alpha$ this low-frequency mode is unstable, and with alpha it is damped less strong than the high-frequency mode. This explains the behaviour discussed in Section 3(b) where the high-frequency mode decays faster than the low-frequency mode.

\end{document}